\documentclass[a4paper]{article}

\usepackage[cp1250]{inputenc}
\usepackage[IL2]{fontenc}
\usepackage{a4wide}
\usepackage[english]{babel}
\usepackage{euscript}
\usepackage{amstext,amsbsy,amscd,amssymb}
\usepackage{amsmath,enumerate}
\usepackage{amsfonts}
\usepackage{graphics}
\usepackage{graphicx}
\usepackage{microtype}
\usepackage{color}
\include{diagxy}
\allowdisplaybreaks

\let\rarr=\rightarrow

\let\veps=\varepsilon

\let\mfrak=\mathfrak
\let\eus=\EuScript

\def\N{\mathbb{N}}
\def\Z{\mathbb{Z}}

\def\C{\mathbb{C}}

\def\End{\mathop {\rm End} \nolimits}
\def\Hom{\mathop {\rm Hom} \nolimits}

\def\ad{\mathop {\rm ad} \nolimits}

\def\GL{\mathop {\rm GL} \nolimits}

\def\id{{\rm id}}

\def\rank{\mathop {\rm rank} \nolimits}

\def\Ind{\mathop {\rm Ind} \nolimits}

\def\htt{\mathop {\rm ht} \nolimits}

\let\g=\gamma
\let\p=\partial

\newsymbol\squares 1003

\long\def\proof #1{\noindent \emph{Proof.}\ #1 \hfill $\squares$

\medskip}

\newcounter{num}[section]
\numberwithin{equation}{section}
\numberwithin{num}{section}

\long\def\definition #1 {\refstepcounter{num} \noindent {\bf
Definition \thenum.} #1

\medskip}

\long\def\theorem #1{\refstepcounter{num} \noindent {\bf Theorem
\thenum.} #1

\medskip}

\long\def\lemma #1{\refstepcounter{num}  \noindent {\bf Lemma
\thenum.} #1

\medskip}

\long\def\proposition #1{\refstepcounter{num}  \noindent {\bf
Proposition \thenum.} #1

\medskip}

\makeatletter
\newcommand*\if@single[3]{%
  \setbox0\hbox{${\mathaccent"0362{#1}}^H$}%
  \setbox2\hbox{${\mathaccent"0362{\kern0pt#1}}^H$}%
  \ifdim\ht0=\ht2 #3\else #2\fi
  }
\newcommand*\rel@kern[1]{\kern#1\dimexpr\macc@kerna}
\newcommand*\widebar[1]{\@ifnextchar^{{\wide@bar{#1}{0}}}{\wide@bar{#1}{1}}}
\newcommand*\wide@bar[2]{\if@single{#1}{\wide@bar@{#1}{#2}{1}}{\wide@bar@{#1}{#2}{2}}}
\newcommand*\wide@bar@[3]{%
  \begingroup
  \def\mathaccent##1##2{%
    \if#32 \let\macc@nucleus\first@char \fi
    \setbox\z@\hbox{$\macc@style{\macc@nucleus}_{}$}%
    \setbox\tw@\hbox{$\macc@style{\macc@nucleus}{}_{}$}%
    \dimen@\wd\tw@
    \advance\dimen@-\wd\z@
    \divide\dimen@ 3
    \@tempdima\wd\tw@
    \advance\@tempdima-\scriptspace
    \divide\@tempdima 10
    \advance\dimen@-\@tempdima
    \ifdim\dimen@>\z@ \dimen@0pt\fi
    \rel@kern{0.6}\kern-\dimen@
    \if#31
      \overline{\rel@kern{-0.6}\kern\dimen@\macc@nucleus\rel@kern{0.4}\kern\dimen@}%
      \advance\dimen@0.4\dimexpr\macc@kerna
      \let\final@kern#2%
      \ifdim\dimen@<\z@ \let\final@kern1\fi
      \if\final@kern1 \kern-\dimen@\fi
    \else
      \overline{\rel@kern{-0.6}\kern\dimen@#1}%
    \fi
  }%
  \macc@depth\@ne
  \let\math@bgroup\@empty \let\math@egroup\macc@set@skewchar
  \mathsurround\z@ \frozen@everymath{\mathgroup\macc@group\relax}%
  \macc@set@skewchar\relax
  \let\mathaccentV\macc@nested@a
  \if#31
    \macc@nested@a\relax111{#1}%
  \else
    \def\gobble@till@marker##1\endmarker{}%
    \futurelet\first@char\gobble@till@marker#1\endmarker
    \ifcat\noexpand\first@char A\else
      \def\first@char{}%
    \fi
    \macc@nested@a\relax111{\first@char}%
  \fi
  \endgroup
}
\makeatother

\newcommand\rsmraise[1]{%
  \ifx#1\displaystyle .8\else
    \ifx#1\textstyle .8\else
      \ifx#1\scriptstyle .6\else
        .45%
      \fi
    \fi
  \fi}

\title{Generalized Verma modules over $U_q(\mfrak{sl}_n(\C))$}

\author{Vyacheslav Futorny, Libor Křižka, Jian Zhang}

\AtEndDocument{\bigskip{\footnotesize%
  (V.\,Futorny) \textsc{Instituto de Matemática e Estatística, Universidade de Sa\~{o} Paulo, Caixa Postal 66281,\\ Sa\~{o} Paulo, CEP 05315-970, Brasil} \par
  \textit{E-mail address}: \texttt{futorny@ime.usp.br} \par
  \addvspace{\medskipamount}
  (L.\,Křižka) \textsc{Instituto de Matemática e Estatística, Universidade de Sa\~{o} Paulo, Caixa Postal 66281,\\ Sa\~{o} Paulo, CEP 05315-970, Brasil} \par
  \textit{E-mail address}: \texttt{krizka.libor@gmail.com} \par
  \addvspace{\medskipamount}
  (J.\,Zhang) \textsc{Instituto de Matemática e Estatística, Universidade de Sa\~{o} Paulo, Caixa Postal 66281,\\ Sa\~{o} Paulo, CEP 05315-970, Brasil} \par
  \textit{E-mail address}: \texttt{zhang@ime.usp.br} \par
}}

\date{}

\begin{document}

\maketitle

\begin{abstract}
We construct realizations of quantum generalized Verma modules for $U_q(\mathfrak{sl}_n(\C))$  by quantum differential operators. Taking the classical limit $q \rightarrow 1$ provides a realization of classical generalized Verma modules for $\mathfrak{sl}_n(\C)$ by differential operators.
\medskip

\medskip

\medskip

\end{abstract}

\thispagestyle{empty}

\tableofcontents


\section*{Introduction}
\addcontentsline{toc}{section}{Introduction}

Generalized Verma modules for complex simple finite-dimensional Lie algebras play an important role in representation theory of Lie algebras. They were first introduced by Garland and Lepowsky in \cite{Garland-Lepowsky1976}. The theory was further developed by many authors, see \cite{Lepowsky1977a}, \cite{Lepowsky1977}, \cite{Rocha-Caridi1980}, \cite{Fernando1990}, \cite{Boe-Collingwood1990}, \cite{Coleman-Futorny1994}, \cite{Futorny1987}, \cite{Futorny1986}, \cite{Drozd-Ovsienko-Futorny1990}, \cite{Irving-Shelton1988}, \cite{Dimitrov-Mathieu-Penkov2000},  \cite{Mazorchuk2000}, \cite{Mazorchuk-Stroppel2008} and references therein. The generalized Verma modules are a natural generalization of the Verma modules defined in \cite{Verma1966}, they are obtained by the \emph{parabolic} induction for a given choice of a parabolic subalgebra. When a parabolic subalgebra coincides with a Borel subalgebra we obtain the corresponding Verma module. The importance of generalized Verma modules was shown in \cite{Fernando1990}, \cite{Futorny1987}, \cite{Coleman-Futorny1994}, \cite{Dimitrov-Mathieu-Penkov2000} by proving that any \emph{weight} (with respect to a fixed Cartan subalgebra) simple module over a complex simple finite-dimensional Lie algebra $\mathfrak{g}$ is either \emph{cuspidal} or a quotient of a certain generalized Verma module, which in turn is obtained by a parabolic induction from the simple weight module over the Levi factor of the parabolic subalgebra. Let us note that the concept of cuspidality depends whether the weight subspaces have finite or infinite dimension \cite{Fernando1990}, \cite{Futorny1987}. Also, the structure theory of generalized Verma modules differs significantly depending on whether the inducing module over the Levi subalgebra is cuspidal or not. The case of cuspidal inducing modules with finite-dimensional weight spaces was fully settled in \cite{Mazorchuk-Stroppel2008} where it was shown that the block of the category of such modules is equivalent to certain blocks of the category $\mathcal{O}$. On the other hand, the classical construction of generalized Verma modules in \cite{Garland-Lepowsky1976} uses finite-dimensional inducing modules over the Levi subalgebra. Such induced modules have certain universal properties but at the same time they are quotients of the corresponding Verma modules.

It is always useful and important to have a concrete realization of simple modules in terms of differential operators. Such realizations for different representations  of $\mathfrak{sl}_n$ can be obtained, for instance, via the embedding into the Witt algebra $W_{n-1}$ \cite{Shen1986}.
The purpose of the present paper is to study quantum deformations of the  generalized Verma modules and construct realizations of these modules (which are simple generically) by quantum differential operators (Theorem \ref{thm:rho action} and Theorem \ref{thm}).  We note that our construction holds for finite and infinite-dimensional inducing modules over parabolic subalgebras.
Similar realizations can be be constructed for quantum groups of all types. Taking the classical limit $q \rightarrow 1$ provides a realization of classical generalized Verma modules by differential operators.

Throughout the article we use the standard notation $\N$ and $\N_0$ for the set of natural numbers and the set of natural numbers together with zero, respectively.


\section{Quantum Weyl algebras}

For $q \in \C^\times$ satisfying $q \neq \pm 1$ and $v \in \C$, the $q$-number $[v]_q$ is defined by
\begin{align}
  [v]_q = {q^v-q^{-v} \over q-q^{-1}}.
\end{align}
If $n \in \N_0$, then we introduce the $q$-factorial $[n]_q!$ by
\begin{align}
  [n]_q!= \prod_{k=1}^n [k]_q.
\end{align}
The $q$-binomial coefficients are defined by the formula
\begin{align}
  {n \brack k}_q = {[n]_q! \over [k]_q![n-k]_q!},
\end{align}
where $n, k \in \N_0$ and $n \geq k$.
\medskip

Let us consider an associative $\C$-algebra $A$. Let $\sigma \colon A \rarr A$ be a $\C$-algebra automorphism. Then a twisted derivation of $A$ relative to $\sigma$ is a linear mapping $D \colon A \rarr A$ satisfying
\begin{align}
  D(ab) = D(a)\sigma(b) + \sigma^{-1}(a) D(b)
\end{align}
for all $a,b \in A$. An element $a \in A$ induces an inner twisted derivation $\ad_\sigma\! a$ relative to $\sigma$ defined by the formula
\begin{align}
  (\ad_\sigma\! a)(b)= a \sigma(b) - \sigma^{-1}(b)a
\end{align}
for all $a,b \in A$. Let us note that also $D_\sigma = \sigma - \sigma^{-1}$ is a twisted derivation of $A$ relative to $\sigma$.
\medskip

\lemma{Let $D$ be a twisted derivation of $A$ relative to $\sigma$. Then we have
\begin{align}
\begin{aligned}
  \sigma \circ \lambda_a &= \lambda_{\sigma(a)} \circ \sigma, \qquad &    D \circ \lambda_a - \lambda_{\sigma^{-1}(a)} \circ D &= \lambda_{D(a)} \circ \sigma, \\
  \sigma \circ \rho_a &= \rho_{\sigma(a)} \circ \sigma, &
   D \circ \rho_a - \rho_{\sigma(a)} \circ D &= \rho_{D(a)} \circ \sigma^{-1}
\end{aligned}
\end{align}
for all $a \in A$, where $\lambda_a$ and $\rho_a$ denote the left and the right multiplications by $a \in A$, respectively.}

\proof{We have
\begin{gather*}
  (\sigma \circ \lambda_a)(b)= \sigma(ab)= \sigma(a)\sigma(b) = (\lambda_{\sigma(a)} \circ \sigma)(b) \\
  (D \circ \lambda_a)(b)= D(ab) = D(a)\sigma(b)+\sigma^{-1}(a)D(b) = (\lambda_{D(a)}\circ \sigma + \lambda_{\sigma^{-1}(a)} \circ D)(b)
\end{gather*}
for all $a,b \in A$.}

Let $V$ be a finite-dimensional complex vector space and let $\C[V]$ be the $\C$-algebra of polynomial functions on $V$. Further, let $\{x_1,x_2,\dots,x_n\}$ be the linear coordinate functions on $V$ with respect to a basis $\{e_1,e_2,\dots,e_n\}$ of $V$. Then there exists a canonical isomorphism of $\C$-algebras $\C[V]$ and $\C[x_1,x_2,\dots,x_n]$.

Let $q \in \C^\times$ satisfies $q \neq \pm 1$. We define a $\C$-algebra automorphism $\gamma_{q,x_i}$ of $\C[V]$ by
\begin{align}
  \gamma_{q,x_i} = q^{x_i \partial_{x_i}}
\end{align}
and a twisted derivation $\partial_{q,x_i}$ of $\C[V]$ relative to $\gamma_{q,x_i}$ through
\begin{align}
  \partial_{q,x_i} = {1 \over x_i} {q^{x_i\partial_{x_i}} - q^{-x_i\partial_{x_i}} \over q - q^{-1}}
\end{align}
for $i=1,2,\dots,n$.
\medskip

\lemma{Let $q \in \C^\times$ satisfies $q \neq \pm 1$. Further, let $D$ be a twisted derivation of $\C[V]$ relative to $\gamma_{q,x_i}$ for some $i=1,2,\dots,n$. Then we have
\begin{align}
  D =   f_i \partial_{q,x_i},
\end{align}
where $f_i \in \C[V]$.}

\proof{For $j=1,2,\dots,n$ satisfying $j \neq i$, we have
\begin{align*}
  D(x_ix_j) &= D(x_i)\gamma_{q,x_i}(x_j)+\gamma_{q,x_i}^{-1}(x_i)D(x_j) =x_jD(x_i)+q^{-1}x_iD(x_j), \\
  D(x_jx_i) &= D(x_j)\gamma_{q,x_i}(x_i)+\gamma_{q,x_i}^{-1}(x_j)D(x_i) =qx_iD(x_j)+x_jD(x_i),
\end{align*}
which implies that $D(x_j)=0$ for all $j=1,2,\dots,n$ such that $j \neq i$. If we set $f_i=D(x_i)$, then we get
\begin{align*}
  (D-f_i\partial_{q,x_i})(x_j)=0
\end{align*}
for all $j=1,2,\dots,n$, which gives us $D=f_i \partial_{q,x_i}$.}

Let $q \in \C^\times$ satisfies $q \neq \pm 1$. Then based on the previous lemma, we define the quantum Weyl algebra $\eus{A}^q_V$ of the complex vector space $V$ as an associative $\C$-subalgebra of $\End \C[V]$ generated by $x_i$, $\partial_{q,x_i}$ and $\gamma_{q,x_i}^{\pm 1}$ for $i=1,2,\dots,n$. Let us note that the definition of $\eus{A}^q_V$ depends on the choice of a basis $\{e_1,e_2,\dots,e_n\}$ of $V$. Moreover, we have the following nontrivial relations
\begin{align}
 \gamma_{q,x_i} x_i  = q x_i \gamma_{q,x_i}, \qquad \gamma_{q,x_i} \partial_{q,x_i} = q^{-1} \partial_{q,x_i}\gamma_{q,x_i}
\end{align}
and
\begin{align}
  \partial_{q,x_i} x_i - q x_i \partial_{q,x_i}=\gamma_{q,x_i}^{-1}, \qquad   \partial_{q,x_i} x_i - q^{-1} x_i \partial_{q,x_i}=\gamma_{q,x_i}
\end{align}
for $i=1,2,\dots,n$.


\section{Generalized Verma modules}


\subsection{Generalized Verma modules over Lie algebras}

Let us consider a finite-dimensional complex semisimple Lie algebra $\mfrak{g}$. Let $\mathfrak{h}$ be a Cartan subalgebra of $\mfrak{g}$. We denote by $\Delta$ the root system of $\mfrak{g}$ with respect to $\mathfrak{h}$, by $\Delta^{+}$ a positive root system in $\Delta$, and by $\Pi\subset\Delta$ the set of simple roots.

Let $\rank \mfrak{g}=r$ and $\Pi= \{\alpha_1,\alpha_2,\dots,\alpha_r\}$. Then we denote by $\alpha_i^\vee \in \mfrak{h}$ the coroot corresponding to the root $\alpha_i$ and by $\omega_i \in \mfrak{h}^*$ the fundamental weight defined through $\langle \omega_i,\alpha_j^\vee \rangle  = \delta_{ij}$ for all $j=1,2,\dots,r$. We also set
\begin{align}
  Q = \sum_{\alpha \in \Delta^+} \Z \alpha = \bigoplus_{i=1}^r \Z \alpha_i \qquad \text{and} \qquad Q_+= \sum_{\alpha \in \Delta^+} \N_0 \alpha = \bigoplus_{i=1}^r \N_0 \alpha_i
\end{align}
and call $Q$ the root lattice and $Q_+$ the positive root lattice. The Cartan matrix $A=(a_{ij})_{1 \leq i,j \leq r}$ of $\mfrak{g}$ is given by $a_{ij} = \langle \alpha_j, \alpha_i^\vee \rangle$.

Further, we denote by $s_i \in \GL(\mfrak{h}^*)$ the reflection about the hyperplane perpendicular to the root $\alpha_i$. Then we obtain $s_i(\alpha_j) = \alpha_j - a_{ij} \alpha_i$. Let $W_\mfrak{g}$ be the Weyl group of $\mfrak{g}$ generated by $s_i$ for $i=1,2,\dots,r$. Then $W_\mfrak{g}$ is a finite Coxeter group with generators $\{s_1,s_2,\dots,s_r\}$ and the relations
\begin{align}
  (s_is_j)^{m_{ij}} = 1,
\end{align}
where $m_{ii}=1$ and $m_{ij}=2,3,4 \text{ or } 6$ for $a_{ij}a_{ji} = 0,1,2 \text{ or } 3$, respectively, provided $i \neq j$. Together with the Weyl group $W_\mfrak{g}$ it is useful to introduce the (generalized) braid group $B_\mfrak{g}$ of $\mfrak{g}$. It is an infinite group with generators $\{T_1,T_2,\dots,T_r\}$ and the braid relations
\begin{align}
  \underbrace{T_i T_j T_j \dotsm}_{m_{ij}} = \underbrace{T_j T_i T_j \dotsm}_{m_{ji}}
\end{align}
for $i \neq j$, where $m_{ij}=m_{ji}$. Let us note that the Weyl group $W_\mfrak{g}$ is the quotient of $B_\mfrak{g}$ under the further relations $T_i^2=1$ for $i=1,2,\dots,r$. For an element $w \in W_\mfrak{g}$ we introduce the length $\ell(w)$ by
\begin{align}
  \ell(w) = |\Delta^+ \cap w(-\Delta^+)|.
\end{align}
Let us note that the length $\ell(w)$ of $w \in W_\mfrak{g}$ is the smallest nonnegative integer $k \in \N_0$ required for an expression of $w$ into the form
\begin{align}
  w = s_{i_1}s_{i_2} \dotsm s_{i_k}, \label{eq:expression w}
\end{align}
where $i_1,i_2,\dots,i_k \in \{1,2,\dots,r\}$. Such an expression is called a \emph{reduced expression} of $w$ if $k=\ell(w)$. It is well known that there exists a unique element $w_0 \in W_\mfrak{g}$ of the maximal length $\ell(w_0) = |\Delta^+|$ called the \emph{longest element}.
\medskip

The standard Borel subalgebra $\mfrak{b}$ of $\mfrak{g}$ is defined through $\mfrak{b} = \mfrak{h} \oplus \mfrak{n}$ with the nilradical $\mfrak{n}$ and the opposite nilradical $\widebar{\mfrak{n}}$ given by
\begin{align}
  \mfrak{n} = \bigoplus_{\alpha \in \Delta^+} \mfrak{g}_\alpha \qquad \text{and} \qquad \widebar{\mfrak{n}} = \bigoplus_{\alpha \in \Delta^+} \mfrak{g}_{-\alpha}. \label{eq:niradical borel}
\end{align}
Moreover, we have a triangular decomposition
\begin{align}
  \mfrak{g}= \widebar{\mfrak{n}} \oplus \mfrak{h} \oplus \mfrak{n} \label{eq:triangular decomposition borel}
\end{align}
of the Lie algebra $\mfrak{g}$.

Further, let us consider a subset $\Sigma$ of $\Pi$ and denote by $\Delta_\Sigma$ the root subsystem in $\mfrak{h}^*$ generated by $\Sigma$. Then the standard parabolic subalgebra $\mfrak{p}=\mfrak{p}_\Sigma$ of $\mfrak{g}$ associated to $\Sigma$ is defined through $\mfrak{p} = \mfrak{l} \oplus \mfrak{u}$ with the nilradical $\mfrak{u}$ and the opposite nilradical $\widebar{\mfrak{u}}$ given by
\begin{align}
  \mfrak{u} = \bigoplus_{\alpha \in \Delta^+ \setminus \Delta_\Sigma} \mfrak{g}_\alpha \qquad \text{and} \qquad \widebar{\mfrak{u}} = \bigoplus_{\alpha \in \Delta^+ \setminus \Delta_\Sigma} \mfrak{g}_{-\alpha} \label{eq:niradical parabolic}
\end{align}
and with the reductive Levi subalgebra $\mfrak{l}$ defined by
\begin{align}
  \mfrak{l} = \mfrak{h} \oplus \bigoplus_{\alpha \in \Delta_\Sigma} \mfrak{g}_\alpha.
\end{align}
Moreover, we have a triangular decomposition
\begin{align}
  \mfrak{g}= \widebar{\mfrak{u}} \oplus \mfrak{l} \oplus \mfrak{u} \label{eq:triangular decomposition parabolic}
\end{align}
of the Lie algebra $\mfrak{g}$. Furthermore, we define the $\Sigma$-height $\htt_\Sigma(\alpha)$ of $\alpha \in \Delta$ by
\begin{align}
  \htt_\Sigma\!\bigg(\sum_{\alpha \in \Pi} m_\alpha \alpha\bigg)=\sum_{\alpha \in \Pi \setminus \Sigma}  m_\alpha.
\end{align}
This gives us a structure of a $|k|$-graded Lie algebra on $\mfrak{g}$ for some $k \in \N_0$. Let us note that if $\Sigma =\emptyset$ then $\mfrak{p} = \mfrak{b}$ and if $\Sigma = \Pi$ then $\mfrak{p} = \mfrak{g}$.
\medskip

\definition{Let $V$ be a simple $\mfrak{p}$-module satisfying $\mfrak{u}V=0$. Then the \emph{generalized Verma module} $M^\mfrak{g}_\mfrak{p}(V)$ is the induced module
\begin{align}
M^\mfrak{g}_\mfrak{p}(V)= \Ind^{U(\mfrak{g})}_{U(\mfrak{p})}(V) \equiv U(\mfrak{g})\otimes_{U(\mfrak{p})}\! V \simeq U(\widebar{\mfrak{u}}) \otimes_\C\! V,
\end{align}
where the last isomorphism of $U(\widebar{\mfrak{u}})$-modules follows from  Poincar\'e--Birkhoff--Witt theorem.}

If $\mfrak{l}$ is the Cartan subalgebra $\mfrak{h}$, then $\mfrak{p}$ is the Borel subalgebra $\mfrak{b}$. In that case, any simple $\mfrak{p}$-module $V$ is $1$-dimensional and $M^\mfrak{g}_\mfrak{p}(V)$ is the corresponding Verma module. Moreover, if $V$ is a finite-dimensional $\mfrak{p}$-module, then $M^\mfrak{g}_\mfrak{p}(V)$ is a homomorphic image of a certain Verma module $M^\mfrak{g}_\mfrak{b}(W)$, where $W$ is a $1$-dimensional $\mfrak{b}$-module. Let us note that $M^\mfrak{g}_\mfrak{p}(V)$ has a unique simple quotient $L^\mfrak{g}_\mfrak{p}(V)$ and generically $M^\mfrak{g}_\mfrak{p}(V) \simeq L^\mfrak{g}_\mfrak{p}(V)$.


\subsection{Generalized Verma modules over quantum groups}

In this section we describe generalized Verma modules for quantum groups. For more detailed information concerning quantum groups see e.g.\ \cite{Kassel1995}, \cite{Klimyk-Schmudgen1997}, \cite{Chari-Pressley1994}. We use the notation introduced in the previous section.
\medskip

Let $\mfrak{g}$ be a finite-dimensional complex semisimple Lie algebra of rank $r$ together with the set of simple roots $\Pi = \{\alpha_1,\alpha_2,\dots,\alpha_r\}$, the Cartan matrix $A=(a_{ij})_{1 \leq i,j \leq r}$ and $d_i = {1\over 2} (\alpha_i,\alpha_i)$ for $i=1,2,\dots,r$, where $(\cdot\,,\cdot)$ is the inner product on $\mfrak{h}^*$ induced by the Cartan--Killing form on $\mfrak{g}$ and normalized so that $(\alpha,\alpha)=2$ for short roots $\alpha \in \Delta^+$.

Let $q \in \C^\times$ satisfies $q^{d_i} \neq \pm 1$ for $i=1,2,\dots,r$. Then the quantum group $U_q(\mfrak{g})$ is a unital associative $\C$-algebra generated by $e_i$, $f_i$, $k_i$, $k_i^{-1}$ for $i=1,2,\dots,r$ subject to the relations
\begin{align}
  \begin{gathered}
    k_ik_i^{-1} = 1, \qquad k_ik_j = k_jk_i, \qquad k_i^{-1}k_i = 1, \\
    k_ie_jk_i^{-1} = q^{a_{ij}} e_j, \qquad [e_i,f_j] = \delta_{ij}\, {k_i - k_i^{-1} \over q_i -q_i^{-1}}, \qquad k_if_jk_i^{-1} = q^{-a_{ij}} f_j
  \end{gathered}
\end{align}
for $i,j=1,2,\dots,r$ and the quantum Serre relations
\begin{align}
    \sum_{k=0}^{1-a_{ij}} (-1)^k {1-a_{ij} \brack k}_{q_i} e_i^{1-a_{ij}-k}e_je_i^k=0, \qquad
    \sum_{k=0}^{1-a_{ij}} (-1)^k {1-a_{ij} \brack k}_{q_i} f_i^{1-a_{ij}-k}f_jf_i^k=0
\end{align}
for $i,j=1,2,\dots,r$ satisfying $i \neq j$, where $q_i = q^{d_i}$ for $i=1,2,\dots,r$.

There is a unique Hopf algebra structure on the quantum group $U_q(\mfrak{g})$ with the coproduct $\Delta \colon U_q(\mfrak{g}) \rarr U_q(\mfrak{g}) \otimes_\C U_q(\mfrak{g})$, the counit $\veps \colon U_q(\mfrak{g}) \rarr \C$ and the antipode $S \colon U_q(\mfrak{g}) \rarr U_q(\mfrak{g})$ given by
\begin{align} \label{eq:Hopf algebra structure U_q(g)}
  \begin{gathered}
    \Delta(e_i)=e_i \otimes k_i + 1 \otimes e_i, \\
 \veps(e_i)=0,  \\
 S(e_i)=-e_ik_i^{-1},
  \end{gathered} \qquad
  \begin{gathered}
    \Delta(k_i)=k_i \otimes k_i, \\
    \veps(k_i) =1, \\
    S(k_i)=k^{-1}_i,
  \end{gathered} \qquad
  \begin{gathered}
    \Delta(f_i)=f_i \otimes 1 + k_i^{-1} \otimes f_i, \\
    \veps(f_i)=0, \\
    S(f_i)=-k_if_i
  \end{gathered}
\end{align}
for $i=1,2,\dots,r$.

Moreover, there exists a homomorphism of the braid group $B_\mfrak{g}$ into the group of $\C$-algebra automorphisms of $U_q(\mfrak{g})$ determined by
\begin{align}
T_i(e_i) = -k_i^{-1}f_i, \qquad T_i(k_j) = k_jk_i^{-a_{ij}}, \qquad T_i(f_i) = -e_ik_i
\end{align}
for $i,j=1,2,\dots,r$ and
\begin{align}
  \begin{aligned}
    T_i(e_j) &= \sum_{s=0}^{-a_{ij}} (-1)^{s-a_{ij}} q_i^{-s} {e_i^s \over [s]_{q_i}!}\,e_j {e_i^{-a_{ij}-s} \over [-a_{ij}-s]_{q_i}!}, \\
    T_i(f_j) &= \sum_{s=0}^{-a_{ij}} (-1)^{s-a_{ij}} q_i^s {f_i^{-a_{ij}-s} \over [-a_{ij}-s]_{q_i}!}\,f_j {f_i^s \over [s]_{q_i}!}
  \end{aligned}
\end{align}
for $i,j=1,2,\dots,r$ satisfying $i \neq j$.

Let $w_0 \in W_\mfrak{g}$ be the longest element in the Weyl group $W_\mfrak{g}$ with a reduced expression
\begin{align}
  w_0 = s_{i_1}s_{i_2}\dotsm s_{i_n},
\end{align}
where $n=|\Delta^+|$. If we set
\begin{align}
  \beta_k = s_{i_1}s_{i_2}\dotsm s_{i_{k-1}}(\alpha_{i_k})
\end{align}
for $k=1,2,\dots,n$, then the sequence $\beta_1,\beta_2,\dots,\beta_n$ exhausts all positive roots $\Delta^+$ of $\mfrak{g}$. Hence, we define
\begin{align}
  e_{\beta_k} = T_{i_1}T_{i_2} \dotsm T_{i_{k-1}}(e_{i_k}) \qquad \text{and} \qquad f_{\beta_k} = T_{i_1}T_{i_2} \dotsm T_{i_{k-1}}(f_{i_k})
\end{align}
and get elements of $U_q(\mfrak{g})$ called root vectors of $U_q(\mfrak{g})$ corresponding to the roots $\beta_k$ and $-\beta_k$ for $k=1,2,\dots,n$, respectively.
\medskip

Let $U_q(\mfrak{n})$ and $U_q(\widebar{\mfrak{n}})$ be the $\C$-subalgebras of $U_q(\mfrak{g})$ generated by the root vectors $e_i$ for $i=1,2,\dots,r$ and $f_i$ for $i=1,2,\dots,r$, respectively. For the quantum group $U_q(\mfrak{g})$ we have a direct sum decomposition
\begin{align}
U_q(\mfrak{g})=\bigoplus_{\alpha\in Q}U_q^\alpha(\mfrak{g}),
\end{align}
where
\begin{align}
U_q^\alpha(\mfrak{g})=\{u\in U_q(\mfrak{g});\, k_i u k_i^{-1}=q^{\langle \alpha, \alpha_i^\vee \rangle}u\ \text{for $i=1,2,\dots,r$}\}.
\end{align}
Since $U_q^\alpha(\mfrak{g}) U_q^\beta(\mfrak{g}) \subset U^{\alpha+\beta}_q(\mfrak{g})$, the preceding shows that $U_q(\mfrak{g})$ is a $Q$-graded $\C$-algebra. Moreover, this grading induces $Q$-grading on the $\C$-subalgebras $U_q(\mfrak{n})$ and $U_q(\widebar{\mfrak{n}})$ as well.
In particular, we have
\begin{align}
U_q(\mathfrak{n})=\bigoplus_{\alpha\in Q_+}U_q^\alpha(\mathfrak{n}) \qquad \text{and} \qquad
U_q(\widebar{\mathfrak{n}})=\bigoplus_{\alpha\in Q_+}U_q^{-\alpha}(\widebar{\mathfrak{n}}),
\end{align}
where $U_q^\alpha(\mfrak{n}) = U_q^\alpha(\mfrak{g}) \cap U_q(\mfrak{n})$ and $U_q^\alpha(\widebar{\mfrak{n}}) = U_q^\alpha(\mfrak{g}) \cap U_q(\widebar{\mfrak{n}})$ for $\alpha \in Q$.

Further, we denote by $U_q(\mathfrak{h})$ and $U_q(\mfrak{b})$ the $\C$-subalgebras of $U_q(\mathfrak{g})$ generated by the elements $k_i$, $k_i^{-1}$ for $i=1,2,\dots,r$ and $e_i$, $k_i$, $k_i^{-1}$ for $i=1,2,\dots,r$, respectively. Then we have
\begin{align}
  U_q(\mfrak{b}) \simeq  U_q(\mfrak{h}) \otimes_\C U_q(\mfrak{n}).
\end{align}
Moreover, we have a triangular decomposition
\begin{align}
  U_q(\mfrak{g}) \simeq U_q(\widebar{\mfrak{n}}) \otimes_\C U_q(\mfrak{h}) \otimes_\C U_q(\mfrak{n})
\end{align}
of the quantum group $U_q(\mfrak{g})$. Let us note that $U_q(\mfrak{h})$ and $U_q(\mfrak{b})$ are Hopf subalgebras of $U_q(\mfrak{g})$ unlike $U_q(\mfrak{n})$ and $U_q(\widebar{\mfrak{n}})$.

Let $\Sigma$ be a subset of $\Pi$. Then we have the standard parabolic subalgebra $\mfrak{p}$ of $\mfrak{g}$ associated to $\Sigma$ with the nilradical $\mfrak{u}$, the opposite nilradical $\widebar{\mfrak{u}}$ and the Levi subalgebra $\mfrak{l}$.

Let $U_q(\mfrak{u})$ and $U_q(\widebar{\mfrak{u}})$ be the $\C$-subalgebras of $U_q(\mfrak{g})$ generated by the root vectors $e_\alpha$ for $\alpha \in \Delta^+$ satisfying $\htt_\Sigma(\alpha) \neq 0$ and $f_\alpha$ for $\alpha \in \Delta^+$ satisfying $\htt_\Sigma(\alpha) \neq 0$, respectively. Further, we denote by $U_q(\mfrak{l})$ the Levi quantum subgroup of $U_q(\mfrak{g})$ generated by the elements $k_i$, $k_i^{-1}$ for $i=1,2,\dots,r$ and the root vectors $e_i$, $f_i$ for $i=1,2,\dots,r$ such that $\alpha_i \in \Sigma$. Finally, we define the parabolic quantum subgroup $U_q(\mfrak{p})$ of $U_q(\mfrak{g})$ as the $\C$-subalgebra of $U_q(\mfrak{g})$ generated by $e_i$, $k_i$ for $i=1,2,\dots,r$ and $f_i$ for $i=1,2,\dots,r$ such that $\alpha_i \in \Sigma$. Then we have
\begin{align}
  U_q(\mfrak{p}) \simeq  U_q(\mfrak{l}) \otimes_\C U_q(\mfrak{u}).
\end{align}
Moreover, we have a triangular decomposition
\begin{align}
  U_q(\mfrak{g}) \simeq U_q(\widebar{\mfrak{u}}) \otimes_\C U_q(\mfrak{l}) \otimes_\C U_q(\mfrak{u})
\end{align}
of the quantum group $U_q(\mfrak{g})$. Let us note that $U_q(\mfrak{l})$ and $U_q(\mfrak{p})$ are Hopf subalgebras of $U_q(\mfrak{g})$ unlike $U_q(\mfrak{u})$ and $U_q(\widebar{\mfrak{u}})$.
\medskip

\definition{Let $V$ be a simple $U_q(\mathfrak{p})$-module satisfying $U_q(\mfrak{u})V=0$. Then the generalized Verma module $M^\mfrak{g}_{\mfrak{p},q}(V)$ is the induced module
\begin{align}
M^\mfrak{g}_{\mfrak{p},q}(V) = \Ind_{U_q(\mathfrak{p})}^{U_q(\mathfrak{g})}(V)
\equiv U_q(\mathfrak{g})\otimes_{U_q(\mathfrak{p})} \!V \simeq U_q(\widebar{\mathfrak{u}})\otimes_\C\! V,
\end{align}
where the last isomorphism of $U_q(\widebar{\mfrak{u}})$-modules follows from Poincar\'e--Birkhoff--Witt theorem.}

It is well known that certain simple highest weight modules for $U_q(\mfrak{g})$ are true deformations of  simple highest weight modules for $\mfrak{g}$ in the sense of Lusztig \cite{Lusztig1988}, that is these modules have the same character formula and the latter can be obtained by the classical limit via the $\mathbb{A}$-forms of $U_q(\mfrak{g})$. We refer to the paper \cite{Melville1999} where the $\mathbb{A}$-forms technique in quantum deformation was described in details.  Using this method one can easily show that some generalized Verma modules for $U_q(\mfrak{g})$ are true deformations of generalized Verma modules for $\mfrak{g}$.


\section{Representations of the quantum group $U_q(\mfrak{sl}_n(\C))$}


\subsection{The quantum group $U_q(\mfrak{sl}_n(\C))$}

Let us consider the finite-dimensional complex simple Lie algebra $\mfrak{sl}_n(\C)$ of rank $n-1$ together with the set of simple roots $\Pi = \{\alpha_1,\alpha_2,\dots,\alpha_{n-1}\}$ and the Cartan matrix $A=(a_{ij})_{1 \leq i,j \leq n-1}$ given by $a_{ii}=2$, $a_{ij}=-1$ if $|i-j|=1$ and $a_{ij}=0$ if $|i-j|> 1$.

Let $q \in \C^\times$ satisfies $q \neq \pm 1$. Then the quantum group $U_q(\mfrak{sl}_n(\C))$ is a unital associative $\C$-algebra generated by $e_i$, $f_i$, $k_i$, $k_i^{-1}$ for $i=1,2,\dots,n-1$ subject to the relations
\begin{align} \label{eq:quantum group U_q(sl(n,C)) rel I}
  \begin{gathered}
    k_ik_i^{-1} = 1, \qquad k_ik_j = k_jk_i, \qquad k_i^{-1}k_i = 1, \\
    k_ie_jk_i^{-1} = q^{a_{ij}} e_j, \qquad [e_i,f_j] = \delta_{ij}\, {k_i - k_i^{-1} \over q -q^{-1}}, \qquad k_if_jk_i^{-1} = q^{-a_{ij}}f_j
  \end{gathered}
\end{align}
for $i,j=1,2,\dots,n-1$ and the quantum Serre relations
\begin{align} \label{eq:quantum group U_q(sl(n,C)) rel II}
\begin{gathered}
e_{i}^2e_{j}-(q+q^{-1})e_ie_je_i+e_je_{i}^2=0, \\
e_{i}e_j=e_je_i,
\end{gathered}
\quad
\begin{gathered}
f_{i}^2f_{j}-(q+q^{-1})f_if_jf_i+f_jf_{i}^2=0  \\
f_if_j=f_jf_i,
\end{gathered}
\begin{gathered}
  \quad (|i-j|=1), \\
  \quad (|i-j|>1).
\end{gathered}
\end{align}
Moreover, there exists a unique Hopf algebra structure on the quantum group $U_q(\mathfrak{sl}_n(\C))$ with the coproduct $\Delta_1 \colon U_q(\mathfrak{sl}_n(\C)) \rarr U_q(\mathfrak{sl}_n(\C)) \otimes_\C U_q(\mathfrak{sl}_n(\C))$, the counit $\veps_1 \colon U_q(\mathfrak{sl}_n(\C)) \rarr \C$ and the antipode $S_1 \colon U_q(\mathfrak{sl}_n(\C))\rarr U_q(\mathfrak{sl}_n(\C))$ given by
\begin{align} \label{eq:Hopf algebra structure U_q(sl(n,C)) I}
\begin{gathered}
 \Delta_1(e_i)=e_i \otimes k_i + 1 \otimes e_i, \\
 \veps_1(e_i)=0,  \\
 S_1(e_i)=-e_i k_i^{-1},
\end{gathered} \qquad
\begin{gathered}
  \Delta_1(k_i)= k_i  \otimes k_i, \\
  \veps_1(k_i) =1, \\
  S_1(k_i)=k_i^{-1},
\end{gathered} \qquad
\begin{gathered}
  \Delta_1(f_i)=f_i \otimes 1 + k_i^{-1} \otimes f_i, \\
  \veps_1(f_i)=0, \\
  S_1(f_i)=-k_if_i
\end{gathered}
\end{align}
for $i=1,2,\dots,n-1$. Let us note that we can introduce a different unique Hopf algebra structure on $U_q(\mfrak{sl}_n(\C))$ with the coproduct $\Delta_2 \colon U_q(\mathfrak{sl}_n(\C)) \rarr U_q(\mathfrak{sl}_n(\C)) \otimes_\C U_q(\mathfrak{sl}_n(\C))$, the counit $\veps_2 \colon U_q(\mathfrak{sl}_n(\C)) \rarr \C$ and the antipode $S_2 \colon U_q(\mathfrak{sl}_n(\C))\rarr U_q(\mathfrak{sl}_n(\C))$ given by
\begin{align} \label{eq:Hopf algebra structure U_q(sl(n,C)) II}
\begin{gathered}
 \Delta_2(e_i)=e_i \otimes k_i^{-1} + 1 \otimes e_i, \\
 \veps_2(e_i)=0,  \\
 S_2(e_i)=-e_ik_i,
\end{gathered} \qquad
\begin{gathered}
  \Delta_2(k_i)= k_i  \otimes k_i, \\
  \veps_2(k_i) =1, \\
  S_2(k_i)=k_i^{-1},
\end{gathered} \qquad
\begin{gathered}
  \Delta_2(f_i)=f_i \otimes 1 + k_i \otimes f_i, \\
  \veps_2(f_i)=0, \\
  S_2(f_i)=-k_i^{-1}f_i
\end{gathered}
\end{align}
for $i=1,2,\dots,n-1$.

Furthermore, there is a homomorphism of the braid group $B_{\mfrak{sl}_n(\C)}$ into the group of $\C$-algebra automorphisms of $U_q(\mfrak{sl}_n(\C))$ determined by
\begin{align}
T_i(e_i) = -f_ik_i^{-1}, \qquad T_i(k_i) =k_i^{-1}, \qquad T_i(f_i) = -k_i e_i
\end{align}
for $i=1,2,\dots,n-1$ and
\begin{align}
\begin{gathered}
  T_i(e_j) = e_ie_j-qe_je_i, \\
  T_i(e_j) =e_j,
\end{gathered} \quad
\begin{gathered}
  T_i(k_j) = k_ik_j, \\
  T_i(k_j) =k_j,
\end{gathered} \quad
\begin{gathered}
  T_i(f_j) = f_jf_i - q^{-1}f_if_j\\
  T_i(f_j) = f_j
\end{gathered} \quad
\begin{gathered}
  (|i-j|=1), \\
  (|i-j|>1).
\end{gathered}
\end{align}
Let us note that a simple computation shows that
\begin{align}
  T_iT_j(e_i) = e_j \qquad \text{and} \qquad T_iT_j(f_i)=f_j
\end{align}
for $i,j=1,2,\dots, n-1$ such that $|i-j|=1$.

Now, we construct root basis of $U_q(\mfrak{n})$ and $U_q(\widebar{\mfrak{n}})$ by the approach described in the previous section. The longest element $w_0$ in the Weyl group $W_{\mfrak{sl}_n(\C)}$ has a reduced expression
\begin{align}
  w_0 = s_1\dotsm s_{n-1} s_1\dotsm s_{n-2} \dotsm s_1s_2s_1.
\end{align}
If we set
\begin{align}
  w_{i,j} = s_1\dotsm s_{n-1} \dotsm s_1\dotsm s_{n-i+1} s_1 \dotsm s_{j-i-1},
\end{align}
we obtain
\begin{align}
  w_{i,j}(\alpha_{j-i}) = \alpha_i + \alpha_{i+1} + \dotsb + \alpha_{j-1}
\end{align}
for $1 \leq i < j \leq n$. Hence, we denote by
\begin{align}
  E_{i,j} = T_{w_{i,j}}(e_{j-i})\qquad \text{and} \qquad E_{j,i} = T_{w_{i,j}}(f_{j-i})
\end{align}
elements of $U_q(\mfrak{n})$ and $U_q(\widebar{\mfrak{n}})$ for $1 \leq i < j \leq n$, respectively, where $T_{w_{i,j}}$ stands for
\begin{align}
  T_{w_{i,j}} = T_1\dotsm T_{n-1} \dotsm T_1\dotsm T_{n-i+1} T_1 \dotsm T_{j-i-1}.
\end{align}
Furthermore, we define by
\begin{align}
  K_{i,j} = k_i k_{i+1} \dotsm k_{j-1}
\end{align}
elements of $U_q(\mfrak{h})$ for $1 \leq i < j \leq n$.
\medskip

\proposition{We have
\begin{align}
E_{i,i+1}=e_i, \qquad \qquad E_{i+1,i}=f_i
\end{align}
for $i=1,2,\dots,n-1$ and
\begin{align}
\begin{aligned}
E_{i,j}&=E_{i,k}E_{k,j}-qE_{k,j}E_{i,k}    \\
E_{i,j}&=E_{i,k}E_{k,j}-q^{-1}E_{k,j}E_{i,k}
\end{aligned}\quad
\begin{aligned}
  \text{ for } 1 \leq i<k<j \leq n, \\
  \text{ for } n \geq i>k>j \geq 1.
\end{aligned}
\end{align}
}

\proof{Let us assume that $i<j$. For $1 \leq i < k \leq n$ we have
\begin{align*}
  T_1 \dotsm T_k (e_i) = T_1 \dotsm T_i T_{i+1} (e_i) = T_1\dotsm T_{i-1}(e_{i+1}) = e_{i+1},
\end{align*}
which implies $E_{i,i+1} = e_i$ for $i=1,2,\dots,n-1$. Further, we may write
\begin{align*}
  E_{i,j} &= T_{w_{i,j}}(e_{j-i}) = T_{w_{i,j-1}}T_{j-i-1}(e_{j-i}) = T_{w_{i,j-1}}(e_{j-i-1}e_{j-i}-qe_{j-i}e_{j-i-1}) \\
  &= T_{w_{i,j-1}}(e_{j-i-1})T_{w_{i,j-1}}(e_{j-i}) - qT_{w_{i,j-1}}(e_{j-i})T_{w_{i,j-1}}(e_{j-i-1}) \\
  &= E_{i,j-1} E_{j-1,j} - q E_{j-1,j} E_{i,j-1}
\end{align*}
for $j-i>1$. Hence, we proved the statement for $j-i=1$ and $j-i=2$. The rest of the proof is by induction on $j-i$. For $j-i>2$ we have $E_{i,j} = E_{i,j-1}E_{j-1,j}-qE_{j-1,j}E_{i,j-1}$, which together with the induction assumption $E_{i,j-1} = E_{i,k} E_{k,j-1} - q E_{k,j-1} E_{i,k}$ for $1 \leq i < k < j-1 < n$ gives us
\begin{align*}
  E_{i,j} &= (E_{i,k} E_{k,j-1} - q E_{k,j-1} E_{i,k})E_{j-1,j}-qE_{j-1,j}(E_{i,k} E_{k,j-1} - q E_{k,j-1} E_{i,k}) \\
  &= E_{i,k}(E_{k,j-1}E_{j-1,j}- qE_{j-1,j}E_{k,j-1}) - q(E_{k,j-1}E_{j-1,j}- qE_{j-1,j}E_{k,j-1})E_{i,k} \\
  & = E_{i,k}E_{k,j}-qE_{k,j}E_{i,k},
\end{align*}
where we used $E_{i,k}E_{j-1,j}=E_{j-1,j}E_{i,k}$ in the second equality. For $i > j$ the proof goes along the same lines. This finishes the proof.}

Let us note that the root vectors $E_{i,j}$ of $U_q(\mfrak{sl}_n(\C))$ coincide with the elements introduced by Jimbo in \cite{Jimbo1986}. Moreover, these vectors are linearly independent in $U_q(\mfrak{sl}_n(\C))$ and they have analogous properties as the corresponding elements $E_{i,j}$, $i,j=1,2,\dots,n$ and $i \neq j$, in the matrix realization of $\mathfrak{sl}_n(\C)$.
\medskip

\lemma{\label{lem:relations}
We have
\begin{align*}
& E_{i,k}E_{k,j}^m=q^{-m}E_{k,j}^mE_{i,k}+[m]_q E_{k,j}^{m-1}E_{i,j}\ \text{ for } i>k>j, \\
& E_{i,k}^mE_{k,j}= q^{-m}E_{k,j}E_{i,k}^m + [m]_qE_{i,j}E_{i,k}^{m-1}\ \text{ for }  i>k>j, \\
& E_{i,k}^m E_{i,j}=q^{m}E_{i,j} E_{i,k}^m\ \text{ for } i>k>j, \\
& E_{i,j} E_{k,j}^m=q^{m}E_{k,j}^m E_{i,j}\ \text{ for } i>k>j, \\
& E_{i,i+1}E_{i+1,i}^{m}=E_{i+1,i}^{m}E_{i,i+1}+[m]_qE_{i+1,i}^{m-1} {q^{-m+1}K_{i,i+1}-q^{m-1}K_{i,i+1}^{-1} \over q-q^{-1}}, \\
& E_{i,j}E_{k,i}^m=E_{k,i}^mE_{i,j}-q^{m-2}[m]_qE_{k,i}^{m-1}E_{k,j}K_{i,j}^{-1}\
\text{ for } i<j<k, \\
& E_{j,k}E_{k,i}^{m}=  E_{k,i}^{m} E_{j,k}+[m]_q  E_{j,i}E_{k,i}^{m-1}K_{j,k}\
\text{ for } i<j<k, \\
& E_{\ell,i}E_{k,j}=E_{k,j}E_{\ell,i}\ \text{ for } i<j<k<\ell, \\
& E_{\ell,j} E_{k,i} - E_{k,i} E_{\ell,j} = (q-q^{-1}) E_{k,j}E_{\ell,i}\ \text{ for } i<j<k<\ell, \\
& E_{\ell,i}E_{j,k}=E_{j,k}E_{\ell,i}\ \text{ for } i<j<k<\ell
\end{align*}
in the quantum group $U_q(\mfrak{sl}_n(\C))$.}

\proof{All formulas are easy to be verified by induction.}

\lemma{\label{lem:coproduct}
We have
\begin{align}
\Delta_2(E_{j,i})= K_{i,j} \otimes E_{j,i} + E_{j,i}\otimes 1
+ (q-q^{-1}) \sum_{i<k<j} E_{k,i}K_{k,j} \otimes E_{j,k}
\end{align}
for $1 \leq i<j\leq n$.}

\proof{We prove the statement by induction on $j-i$. The case $j-i=1$ follows immediately from \eqref{eq:Hopf algebra structure U_q(sl(n,C)) II}. Further, for $j-i>1$ we have $E_{j,i} =E_{j,i+1}E_{i+1,i}-q^{-1}E_{i+1,i}E_{j,i+1}$.
Therefore, we may write $\Delta_2(E_{j,i})=\Delta_2(E_{j,i+1}) \Delta_2(E_{i+1,i})
-q^{-1} \Delta_2(E_{i+1,i}) \Delta_2(E_{j,i+1})$.
By induction assumption we have
\begin{align*}
\Delta_2(E_{j,i+1})= K_{i+1,j} \otimes E_{j,i+1} + E_{j,i+1} \otimes 1  +(q-q^{-1})\sum_{i+1<k<j} E_{k,i+1}K_{k,j} \otimes E_{j,k}
\end{align*}
and also
\begin{align*}
  \Delta_2(E_{i+1,i})= K_{i,i+1} \otimes E_{i+1,i} + E_{i+1,i} \otimes 1,
\end{align*}
which gives us
\begin{align*}
  \Delta_2(E_{j,i})&= E_{j,i+1}E_{i+1,i} \otimes 1 - q^{-1} E_{i+1,i} E_{j,i+1} \otimes 1 \\
  & \quad + K_{i+1,j}E_{i+1,i} \otimes E_{j,i+1} - q^{-1} E_{i+1,i} K_{i+1,j} \otimes E_{j,i+1} \\
  & \quad + E_{j,i+1} K_{i,i+1} \otimes E_{i+1,i} - q^{-1} K_{i,i+1} E_{j,i+1}  \otimes E_{i+1,i} \\
  & \quad + K_{i+1,j} K_{i,i+1} \otimes E_{j,i+1} E_{i+1,i} - q^{-1} K_{i,i+1} K_{i+1,j} \otimes E_{i+1,i} E_{j,i+1} \\
  & \quad + (q-q^{-1})\sum_{i+1<k<j} E_{k,i+1}K_{k,j} E_{i+1,i} \otimes  E_{j,k} \\
  & \quad - q^{-1} (q-q^{-1})\sum_{i+1<k<j} E_{i+1,i} E_{k,i+1}K_{k,j} \otimes E_{j,k} \\
  & \quad + (q-q^{-1})\sum_{i+1<k<j} E_{k,i+1}K_{k,j} K_{i,i+1} \otimes  E_{j,k} E_{i+1,i} \\
  &\quad- q^{-1} (q-q^{-1})\sum_{i+1<k<j} K_{i,i+1} E_{k,i+1}K_{k,j} \otimes E_{i+1,i} E_{j,k}.
\end{align*}
Further, using the relations $K_{i+1,j}E_{i+1,i}=qE_{i+1,i}K_{i+1,j}$, $K_{i,i+1}E_{j,i+1}=qE_{j,i+1}K_{i,i+1}$ and $K_{i,i+1}E_{k,i+1}=qE_{k,i+1}K_{i,i+1}$ we may simplified $\Delta_2(E_{j,i})$ into the form
\begin{align*}
  \Delta_2(E_{j,i})&= (E_{j,i+1}E_{i+1,i} - q^{-1} E_{i+1,i} E_{j,i+1}) \otimes 1
   + K_{i,j} \otimes (E_{j,i+1} E_{i+1,i} - q^{-1} E_{i+1,i} E_{j,i+1}) \\
  & \quad + (q- q^{-1}) E_{i+1,i} K_{i+1,j} \otimes E_{j,i+1} \\
  & \quad + (q-q^{-1})\sum_{i+1<k<j} (E_{k,i+1}E_{i+1,i}-q^{-1} E_{i+1,i} E_{k,i+1}) K_{k,j}\otimes  E_{j,k}.
\end{align*}
Therefore, we have
\begin{align*}
 \Delta_2(E_{j,i}) = K_{i,j} \otimes E_{j,i}
  + E_{j,i} \otimes 1 + (q-q^{-1}) \sum_{i<k<j} E_{k,i}K_{k,j} \otimes E_{j,k},
\end{align*}
which finishes the proof.}

\vspace{-2mm}


\subsection{The parabolic induction for $U_q(\mathfrak{sl}_{n+m}(\C))$}

For simplicity we concentrate now on one particular choice of a parabolic quantum subgroup of $U_q(\mathfrak{sl}_{n+m}(\C))$. This offers a good insight into the construction for a general case.

Let $\Sigma=\{\alpha_1,\dots, \alpha_{n-1},\alpha_{n+1},\dots, \alpha_{n+m-1}\}$ be a subset of $\Pi=\{\alpha_1,\alpha_2,\dots,\alpha_{n+m-1}\}$ and let $\mathfrak{p}=\mathfrak{l} \oplus \mathfrak{u}$ be the corresponding parabolic subalgebra of $\mfrak{g}=\mathfrak{sl}_{n+m}(\C)$ with the nilradical $\mfrak{u}$, the opposite nilradical $\widebar{\mfrak{u}}$ and the Levi subalgebra $\mfrak{l}$. We have a triangular decomposition
\begin{align}
  \mfrak{g} = \widebar{\mfrak{u}} \oplus \mfrak{l} \oplus \mfrak{u}
\end{align}
of the Lie algebra $\mfrak{g}$, where $\mfrak{l} \simeq \mfrak{sl}_n(\C) \oplus \mfrak{sl}_m(\C) \oplus \C$, $\widebar{\mfrak{u}} \simeq \Hom(\C^m,\C^n)$ and $\mfrak{u} \simeq \Hom(\C^n,\C^m)$. Furthermore, we have the corresponding quantum parabolic subgroup $U_q(\mfrak{p})$ of $U_q(\mfrak{g})$ with the $\C$-subalgebras $U_q(\mfrak{u})$, $U_q(\widebar{\mfrak{u}})$ and the Levi quantum subgroup $U_q(\mfrak{l})$. Moreover, we have a triangular decomposition
\begin{align}
  U_q(\mfrak{g}) \simeq U_q(\widebar{\mfrak{u}}) \otimes_\C U_q(\mfrak{l}) \otimes_\C U_q(\mfrak{u})
\end{align}
of the quantum group $U_q(\mfrak{g})$.
\medskip

Let $V$ be a $U_q(\mfrak{p})$-module. Then for the induced module $M^\mfrak{g}_{\mfrak{p},q}(V)$ we have
\begin{align}
  U_q(\mfrak{g}) \otimes_{U_q(\mfrak{p})}\! V \simeq U_q(\widebar{\mfrak{u}}) \otimes_\C\! V,
\end{align}
where the isomorphism of vector spaces is in fact an isomorphism of $U_q(\widebar{\mfrak{u}})$-modules. Hence, the action of $U_q(\widebar{\mfrak{u}})$ on $U_q(\widebar{\mfrak{u}}) \otimes_\C \! V$ is just the left multiplication, like in the classical case. Our next step is to describe the action of the Levi quantum subgroup $U_q(\mfrak{l})$ on $U_q(\widebar{\mfrak{u}}) \otimes_\C \! V$, since in the classical case the action of the Levi subalgebra $\mfrak{l}$ on $U(\widebar{\mfrak{u}}) \otimes_\C \! V$ is equal to the tensor product of the adjoint action on $U(\widebar{\mfrak{u}})$ and the action on $V$.
\medskip

Let us recall that the Levi quantum subgroup $U_q(\mfrak{l})$ of $U_q(\mfrak{g})$ has a Hopf algebra structure determined either by \eqref{eq:Hopf algebra structure U_q(sl(n,C)) I} or by \eqref{eq:Hopf algebra structure U_q(sl(n,C)) II}. However, we introduce a different (mixed) Hopf algebra structure on $U_q(\mfrak{l})$ with the coproduct $\Delta \colon U_q(\mfrak{l}) \rarr U_q(\mfrak{l}) \otimes_\C U_q(\mfrak{l})$, the counit $\veps \colon U_q(\mfrak{l}) \rarr U_q(\mfrak{l})$ and the antipode $S \colon U_q(\mfrak{l}) \rarr U_q(\mfrak{l})$ given by
\begin{align} \label{eq:Hopf algebra structure U_q(l) a}
\begin{gathered}
 \Delta(e_i)=e_i \otimes k_i^{-1} + 1 \otimes e_i, \\
 \veps(e_i)=0,  \\
 S(e_i)=-e_ik_i,
\end{gathered} \qquad
\begin{gathered}
  \Delta(k_i)= k_i  \otimes k_i, \\
  \veps(k_i) =1, \\
  S(k_i)=k_i^{-1},
\end{gathered} \qquad
\begin{gathered}
  \Delta(f_i)=f_i \otimes 1 + k_i \otimes f_i, \\
  \veps(f_i)=0, \\
  S(f_i)=-k_i^{-1}f_i
\end{gathered}
\end{align}
for $i=1,2,\dots,n-1$,
\begin{align} \label{eq:Hopf algebra structure U_q(l) c}
\begin{gathered}
  \Delta(k_n) = k_n \otimes k_n, \\
  \veps(k_n) = 1, \\
  S(k_n) = k_n^{-1},
\end{gathered}
\end{align}
and
\begin{align} \label{eq:Hopf algebra structure U_q(l) b}
\begin{gathered}
 \Delta(e_i)=e_i \otimes k_i + 1 \otimes e_i, \\
 \veps(e_i)=0,  \\
 S(e_i)=-e_i k_i^{-1},
\end{gathered} \qquad
\begin{gathered}
  \Delta(k_i)= k_i  \otimes k_i, \\
  \veps(k_i) =1, \\
  S(k_i)=k_i^{-1},
\end{gathered} \qquad
\begin{gathered}
  \Delta(f_i)=f_i \otimes 1 + k_i^{-1} \otimes f_i, \\
  \veps(f_i)=0, \\
  S(f_i)=-k_if_i
\end{gathered}
\end{align}
for $i=n+1,n+2,\dots,n+m-1$.

The Hopf algebra structure on $U_q(\mfrak{l})$ ensures that we can define the (left) adjoint action of $U_q(\mfrak{l})$ on $U_q(\mfrak{g})$ by
\begin{align}
  \ad(a)b = \sum a_{(1)}b S(a_{(2)}), \label{eq:adjoint action}
\end{align}
where
\begin{align}
  \Delta(a) = \sum a_{(1)} \otimes a_{(2)},
\end{align}
for all $a \in U_q(\mfrak{l})$ and $b \in U_q(\mfrak{g})$. Let us note that we also have
\begin{align}
  \ad(a)bc= \sum (\ad(a_{(1)}) b) (\ad(a_{(2)})c)  \label{eq:adjoint action product}
\end{align}
for all $a \in U_q(\mfrak{l})$ and $b,c \in U_q(\mfrak{g})$.
\medskip

\proposition{\label{prop:adjoint action}
The $\C$-subalgebra $U_q(\widebar{\mfrak{u}})$ of $U_q(\mfrak{g})$ is a $U_q(\mfrak{l})$-module with respect to the adjoint action. Moreover, we have
\begin{align}\label{eq:adjoin action levi U_q(l)_a}
\begin{gathered}
  \ad(e_i)E_{n+j,k} = -q^{-1}\delta_{i,k}E_{n+j,k+1}, \qquad  \ad(f_i)E_{j+n,k}= -q\delta_{i+1,k} E_{n+j,k-1}, \\
  \ad(k_i)E_{n+j,k} = q^{-\delta_{i,k}+\delta_{i+1,k}}E_{n+j,k}
\end{gathered}
\end{align}
for $i=1,2,\dots,n-1$,
\begin{align}\label{eq:adjoin action levi U_q(l)_c}
  \ad(k_n)E_{n+j,k} = q^{-\delta_{1,j}-\delta_{n,k}}E_{n+j,k}
\end{align}
and
\begin{align}\label{eq:adjoin action levi U_q(l)_b}
  \begin{gathered}
    \ad(e_{n+i})E_{n+j,k}=\delta_{i+1,j} E_{n+j-1,k}, \qquad \ad(f_{n+i})E_{n+j,k}=\delta_{i,j} E_{n+j+1,k}, \\
    \ad(k_{n+i})E_{n+j,k}= q^{\delta_{i,j}-\delta_{i+1,j}}E_{n+j,k}
  \end{gathered}
\end{align}
for $i=1,2,\dots,m-1$, where $1 \leq j \leq m$ and $1 \leq k \leq n$.}

\proof{Due to the formula \eqref{eq:adjoint action product}, it is enough to verify that the set of generators $\{E_{n+j,k};\, 1 \leq j \leq m, 1 \leq k \leq n\}$ of the $\C$-subalgebra $U_q(\widebar{\mfrak{u}})$ of $U_q(\mfrak{g})$ is preserved by $U_q(\mfrak{l})$ with respect to the adjoint action. The formulas \eqref{eq:adjoin action levi U_q(l)_a}, \eqref{eq:adjoin action levi U_q(l)_b} and \eqref{eq:adjoin action levi U_q(l)_c} are easy consequence of Lemma \ref{lem:relations}.}

\proposition{\label{prop:levi action}
Let $V$ be a $U_q(\mfrak{p})$-module. Then the $U_q(\mfrak{l})$-module structure on $M^\mfrak{g}_{\mfrak{p},q}(V)$ is given by
\begin{align}
  a (u \otimes v) = \sum (\ad a_{(1)})u \otimes a_{(2)}v, \label{eq:action levi}
\end{align}
where
\begin{align}
  \Delta(a) = \sum a_{(1)} \otimes a_{(2)},
\end{align}
for $a \in U_q(\mfrak{l})$, $u \in U_q(\widebar{\mfrak{u}})$ and $v \in V$. In particular, we get that $M^\mfrak{g}_{\mfrak{p},q}(V)$ is isomorphic to $U_q(\widebar{\mfrak{u}}) \otimes_\C \!V$ as $U_q(\mfrak{l})$-module, where the $U_q(\mfrak{l})$-module structure on $U_q(\widebar{\mfrak{u}})$ is given through the adjoint action.}

\proof{For an element $a \in U_q(\mfrak{l})$ we have $\Delta(a) = \sum a_{(1)} \otimes a_{(2)}$, $\Delta(a_{(1)}) = \sum a_{(11)} \otimes a_{(12)}$ and $\Delta(a_{(2)}) = \sum a_{(21)} \otimes a_{(22)}$. Then for $u \in U_q(\widebar{\mfrak{u}})$ and $v \in V$ we may write
\begin{align*}
  \sum (\ad a_{(1)})u \otimes a_{(2)}v &= \sum a_{(11)}uS(a_{(12)}) \otimes a_{(2)}v = \sum a_{(11)}uS(a_{(12)})a_{(2)} \otimes v \\
  & = \sum a_{(1)}uS(a_{(21)})a_{(22)} \otimes v = \sum a_{(1)} u \veps(a_{(2)}) \otimes v  \\
  & = \sum a_{(1)}\veps(a_{(2)}) u  \otimes v = a u \otimes v,
\end{align*}
where we used $(\Delta \otimes \id) \circ \Delta = (\id \otimes \Delta) \circ \Delta$ in the third equality, $m \circ (S \otimes \id) \circ \Delta = i \circ \veps$ in the fourth equality, and $(\id \otimes \veps) \circ \Delta = \id$ in the last equality. Since $U_q(\widebar{\mfrak{u}})$ is a $U_q(\mfrak{l})$-module by Proposition \ref{prop:adjoint action}, we immediately obtain that $M^\mfrak{g}_{\mfrak{p},q}(V)$ is isomorphic to $U_q(\widebar{\mfrak{u}}) \otimes_\C \! V$ as $U_q(\mfrak{l})$-module.}

Let us note that the formula \eqref{eq:action levi} holds for an arbitrary Hopf algebra structure on $U_q(\mfrak{l})$. However, the main difficulty is to find such a Hopf algebra structure that $U_q(\widebar{\mfrak{u}})$ is a $U_q(\mfrak{l})$-module with respect to the adjoint action \eqref{eq:adjoint action}.
\medskip

As a consequence of Proposition \ref{prop:adjoint action} we have that the vector space
\begin{align}
\widebar{\mfrak{u}}_q= \langle \{E_{n+j,k};\, 1 \leq j \leq m, 1 \leq k \leq n\}\rangle
\end{align}
is a $U_q(\mfrak{l})$-submodule of $U_q(\widebar{\mfrak{u}})$. By the specialization $q \rarr 1$ of the root vectors $E_{n+j,k}$, we obtain the canonical root vectors $x_{j,k}$ of $\widebar{\mfrak{u}}$ for $1 \leq j \leq m$ and $1 \leq k \leq n$. Hence, we define an isomorphism $\psi_q \colon \widebar{\mfrak{u}} \rarr \widebar{\mfrak{u}}_q$ of vector spaces by
\begin{align}
  \psi_q(x_{j,k}) = E_{n+j,k} \label{eq:isomorphism u and u_q}
\end{align}
for $1 \leq j \leq m$ and $1 \leq k \leq n$. Let us note that $x=(x_{j,k})_{1 \leq j \leq m, 1 \leq k \leq n}$ gives us linear coordinate functions on $\widebar{\mfrak{u}}^*$.
Further, we introduce a $U_q(\mfrak{l})$-module structure on $\widebar{\mfrak{u}}$ through $\tau_q \colon U_q(\mfrak{l}) \rarr \End \widebar{\mfrak{u}}$ defined by
\begin{align}
  \tau_q(a) = \psi_q^{-1} \circ \ad(a) \circ \psi_q
\end{align}
for all $a \in U_q(\mfrak{l})$. Moreover, when $q$ is specialized to $1$, we get the original $\mfrak{l}$-module structure on $\widebar{\mfrak{u}}$.

For now, let us assume that $q$ is not a root of unity. Then we have $\widebar{\mfrak{u}} \simeq \smash{L^\mfrak{l}_{\mfrak{b} \cap \mfrak{l},q}}(\omega_{n-1}-2\omega_n+ \omega_{n+1})$ as $U_q(\mfrak{l})$-modules, where $\smash{L^\mfrak{l}_{\mfrak{b} \cap \mfrak{l},q}}(\lambda)$ is the simple highest weight $U_q(\mfrak{l})$-module with highest weight $q^\lambda$ for $\lambda \in \mfrak{h}^*$. Further, since we have
\begin{align}
\begin{aligned}
\widebar{\mfrak{u}} \otimes_\C \widebar{\mfrak{u}} &\simeq L^\mfrak{l}_{\mfrak{b} \cap \mfrak{l},q}(2\omega_{n-1}-4\omega_n + 2\omega_{n+1}) \oplus L^\mfrak{l}_{\mfrak{b} \cap \mfrak{l},q}(\omega_{n-2}-2\omega_n +\omega_{n+2}) \\
& \quad \oplus L^\mfrak{l}_{\mfrak{b} \cap \mfrak{l},q}(2\omega_{n-1} -3\omega_n + \omega_{n+2}) \oplus L^\mfrak{l}_{\mfrak{b} \cap \mfrak{l},q}(\omega_{n-2}-3\omega_n + 2 \omega_{n+1})
\end{aligned}
\end{align}
as $U_q(\mfrak{l})$-modules, we define
\begin{align}
  S_q(\widebar{\mfrak{u}}) = T(\widebar{\mfrak{u}})/I_q,
\end{align}
where $I_q$ is the two-sided ideal of the tensor algebra $T(\widebar{\mfrak{u}})$ generated by 
\begin{align}
  L^\mfrak{l}_{\mfrak{b} \cap \mfrak{l},q}(\omega_{n-2}-3\omega_n + 2 \omega_{n+1}) = \langle v_{i,k,\ell}^+, w_{i,j,k,\ell}^+;\, 1\leq i < j \leq m, 1 \leq k < \ell \leq n \rangle,
\end{align}
where
\begin{align}
\begin{aligned}
  v_{i,k,\ell}^+ &= x_{i,\ell}\otimes x_{i,k}-qx_{i,k} \otimes x_{i,\ell}, \\
  w_{i,j,k\ell}^+ &= x_{j,\ell} \otimes x_{i,k}-x_{i,k} \otimes x_{j,\ell}-q x_{j,k} \otimes x_{i,\ell} + q^{-1} x_{i,\ell} \otimes x_{j,k},
\end{aligned}
\end{align}
and by
\begin{align}
  L^\mfrak{l}_{\mfrak{b} \cap \mfrak{l},q}(2\omega_{n-1} -3\omega_n + \omega_{n+2}) = \langle v_{i,j,k}^-, w_{i,j,k,\ell}^-;\, 1\leq i < j \leq m, 1 \leq k < \ell \leq n \rangle,
\end{align}
where
\begin{align}
  \begin{aligned}
    v_{i,j,k}^- & = x_{j,k}\otimes x_{i,k}-qx_{i,k} \otimes x_{j,k}, \\
    w_{i,j,k,\ell}^- & = x_{j,\ell} \otimes x_{i,k} - x_{i,k} \otimes x_{j,\ell} +q^{-1} x_{j,k} \otimes x_{i,\ell} - q x_{i,\ell} \otimes x_{j,k},
  \end{aligned}
\end{align}
which gives us
\begin{align}
  S_q(\widebar{\mfrak{u}}) \simeq \C_q[\widebar{\mfrak{u}}^*]
\end{align}
with
\begin{align}
\begin{aligned}
 \C_q[\widebar{\mfrak{u}}^*]= \C\langle x \rangle /&(x_{i,\ell}x_{i,k}-qx_{i,k}x_{i,\ell},x_{j,k}x_{i,k}-qx_{i,k}x_{j,k}, x_{j,k}x_{i,\ell}-x_{i,\ell}x_{j,k}, \\
 &\  x_{j,\ell}x_{i,k}-x_{i,k}x_{j,\ell}-(q-q^{-1})x_{i,\ell}x_{j,k};\, 1\leq i<j \leq m,1 \leq k < \ell \leq n).
\end{aligned}
\end{align}
In the previous discussion, we assumed that $q$ is not a root of unity. However, the definition of $S_q(\widebar{\mfrak{u}})$ makes sense for all $q \in \C^\times$ satisfying $q \neq \pm 1$.
Moreover, since the two-sided ideal $I_q$ is a $U_q(\mfrak{l})$-submodule of $T(\widebar{\mfrak{u}})$, we obtain that also $S_q(\widebar{\mfrak{u}})$ is a $U_q(\mfrak{l})$-module for all $q \in \C^\times$ satisfying $q \neq \pm 1$. The specialization $q \rarr 1$ gives us $I_q \rarr I$, hence $\C_q[\widebar{\mfrak{u}}^*] \rarr \C[\widebar{\mfrak{u}}^*]$. Let us note that the $\C$-algebra $\C_q[\widebar{\mfrak{u}}^*]$ is usually called the coordinate algebra of the quantum vector space $\widebar{\mfrak{u}}^*$ introduced in \cite{Reshetikhin-Takhtajan-Faddeev1990}.

It follows immediately from Lemma \ref{lem:relations} that the mapping \eqref{eq:isomorphism u and u_q} may be uniquely extended to a $\C$-algebra homomorphism
\begin{align}
  \psi_q \colon \C_q[\widebar{\mfrak{u}}^*] \rarr U_q(\widebar{\mfrak{u}}). \label{eq:symmetrization}
\end{align}
Moreover, since the set
\begin{align}
\{E_{n+1,1}^{r_{1,1}}E_{n+1,2}^{r_{1,2}}\dotsm E_{n+1,n}^{r_{1,n}}\dotsm E_{n+m,1}^{r_{m,1}} \dotsm E_{n+m,n}^{r_{m,n}};\, r \in M_{m,n}(\N_0)\}
\end{align}
forms a basis of $U_q(\widebar{\mathfrak{u}})$, we obtain that $\psi_q$ is an isomorphism of $\C$-algebras. Further, by the formula \eqref{eq:adjoint action product} and the fact that $\psi_q \colon \C_q[\widebar{\mfrak{u}}^*] \rarr U_q(\widebar{\mfrak{u}})$ is an isomorphism of  $\C$-algebras, we get that $\psi_q$ is an isomorphism of $U_q(\mfrak{l})$-modules.

For an $(m \times n)$-matrix $r=(r_{i,j})_{1 \leq i \leq m,1\leq j \leq n}$ with nonnegative integer entries we denote by $x^r$ an element of $\C_q[\widebar{\mfrak{u}}^*]$ defined by
\begin{align}
  x^r = x_{1,1}^{r_{1,1}}x_{1,2}^{r_{1,2}}\dotsm x_{1,n}^{r_{1,n}}\dotsm x_{m,1}^{r_{m,1}} \dotsm x_{m,n}^{r_{m,n}}
\end{align}
and by $E^r$ an element of $U_q(\widebar{\mfrak{u}})$ defined by
\begin{align*}
  E^r = E_{n+1,1}^{r_{1,1}} E_{n+1,2}^{r_{1,2}} \dotsm E_{n+1,n}^{r_{1,n}} \dotsm E_{n+m,1}^{r_{m,1}} \dotsm E_{n+m,n}^{r_{m,n}}.
\end{align*}
Since the $\C$-algebra $\C_q[\widebar{\mfrak{u}}^*]$ has a basis $\{x^r;\, r \in M_{m,n}(\N_0)\}$ we can find a family of isomorphisms $\varphi_q \colon \C[\widebar{\mfrak{u}}^*] \rarr \C_q[\widebar{\mfrak{u}}^*]$ of vector spaces such that $\varphi_q \rarr \id$ for $q \rarr 1$. Let us define $\varphi_q \colon \C[\widebar{\mfrak{u}}^*] \rarr \C_q[\widebar{\mfrak{u}}^*]$ by
\begin{align}
  \varphi_q(x^r) = x^r
\end{align}
for all $r \in M_{m,n}(\N_0)$. Furthermore, we denote by $1_{i,j} \in M_{m,n}(\N_0)$ the $(m \times n)$-matrix having $1$ at the intersection of the $i$-th row and $j$-th column and $0$ elsewhere. Then the corresponding $U_q(\mfrak{l})$-module structure on $\C[\widebar{\mfrak{u}}^*]$ is given through the homomorphism
\begin{align}
  \rho_q   \colon U_q(\mfrak{l}) \rarr \eus{A}^q_{\widebar{\mfrak{u}}^*}
\end{align}
of associative $\C$-algebras, where $\eus{A}^q_{\widebar{\mfrak{u}}^*}$ is the quantum Weyl algebra of the vector space $\widebar{\mfrak{u}}^*$, defined by
\begin{align}
   \rho_q(a)  = \varphi_q^{-1} \circ \tau_q(a) \circ \varphi_q \label{eq:rho action}
\end{align}
for all $a \in U_q(\mfrak{l})$.
\medskip

Let $\gamma_{i,j}$ be the $\C$-algebra automorphism of $\C[\widebar{\mfrak{u}}^*]$ given by
$\gamma_{i,j} \colon x_{k,\ell} \mapsto q^{\delta_{ik}\delta_{j\ell}} x_{k,\ell}$ and $\partial_{i,j}$ the corresponding twisted derivation of $\C[\widebar{\mfrak{u}}^*]$ relative to $\gamma_{i,j}$ for $1 \leq i \leq m$ and $1 \leq j \leq n$.
\medskip

\theorem{\label{thm:rho action}
We have
\begin{align}\label{eq:rho action U_q(sl(n))}
\begin{gathered}
  \rho_q(e_i) = - \sum_{k=1}^m \prod_{t=k}^m \gamma_{t,i} \gamma_{t,i+1}^{-1} x_{k,i+1} \partial_{k,i}, \qquad \rho_q(f_i) = -\sum_{k=1}^m x_{k,i} \partial_{k,i+1} \prod_{t=1}^k \gamma_{t,i}^{-1} \gamma_{t,i+1}, \\
  \rho_q(k_i) = \prod_{t=1}^m \gamma_{t,i}^{-1}\gamma_{t,i+1}
\end{gathered}
\end{align}
for $i=1,2,\dots,n-1$,
\begin{align}
  \rho_q(k_n) =  \prod_{t=1}^n \gamma_{1,t}^{-1}  \prod_{s=1}^m \gamma_{s,n}^{-1} \label{eq:rho action k_n}
\end{align}
and
\begin{align}\label{eq:rho action U_q(sl(m))}
\begin{gathered}
  \rho_q(e_{n+i}) = \sum_{k=1}^n x_{i,k} \partial_{i+1,k} \prod_{t=k+1}^n \gamma_{i,t} \gamma_{i+1,t}^{-1}, \qquad \rho_q(f_{n+i}) = \sum_{k=1}^n \prod_{t=1}^{k-1} \gamma_{i,t}^{-1}\gamma_{i+1,t} x_{i+1,k} \partial_{i,k}, \\
  \rho_q(k_{n+i}) = \prod_{t=1}^n \gamma_{i,t}\gamma_{i+1,t}^{-1}
\end{gathered}
\end{align}
for $i=1,2,\dots,m-1$.}

\proof{The proof is a straightforward computation. Using \eqref{eq:adjoin action levi U_q(l)_a} and \eqref{eq:Hopf algebra structure U_q(l) a} we have
\begin{align*}
  \tau_q(e_i)(x^r) &= -q^{-1}\sum_{k=1}^m \sum_{j=1}^{r_{k,i}} q^{\sum_{t = k}^m (r_{t,i}-r_{t,i+1})-j}\, x_{1,1}^{r_{1,1}}\dotsm x_{k,i}^{j-1}x_{k,i+1} x_{k,i}^{r_{k,i}-j}\dotsm x_{m,n}^{r_{m,n}} \\
  &= -q^{-1}\sum_{k=1}^m \sum_{j=1}^{r_{k,i}} q^{\sum_{t = k}^m (r_{t,i}-r_{t,i+1}) + r_{k,i}-2j}\, x_{1,1}^{r_{1,1}}\dotsm x_{k,i}^{r_{k,i}-1}x_{k,i+1}^{r_{k,i+1}+1}\dotsm x_{m,n}^{r_{m,n}} \\
  &= -\sum_{k=1}^m q^{\sum_{t = k}^m (r_{t,i}-r_{t,i+1}) -2}[r_{k,i}]_q\, x_{1,1}^{r_{1,1}}\dotsm x_{k,i}^{r_{k,i}-1}x_{k,i+1}^{r_{k,i+1}+1}\dotsm x_{m,n}^{r_{m,n}},
\end{align*}
\begin{align*}
  \tau_q(f_i)(x^r) &= -q \sum_{k=1}^m \sum_{j=1}^{r_{k,i+1}} q^{\sum_{t=1}^k (r_{t,i+1}-r_{t,i})-j} \, x_{1,1}^{r_{1,1}}\dotsm x_{k,i+1}^{r_{k,i+1}-j}x_{k,i} x_{k,i+1}^{j-1}\dotsm x_{m,n}^{r_{m,n}} \\
  &= -q \sum_{k=1}^m \sum_{j=1}^{r_{k,i+1}} q^{\sum_{t=1}^k (r_{t,i+1}-r_{t,i}) + r_{k,i+1}-2j} \, x_{1,1}^{r_{1,1}}\dotsm x_{k,i}^{r_{k,i}+1}x_{k,i+1}^{r_{k,i+1}-1}\dotsm x_{m,n}^{r_{m,n}} \\
  &= - \sum_{k=1}^m  q^{\sum_{t=1}^k (r_{t,i+1}-r_{t,i})}[r_{k,i+1}]_q \, x_{1,1}^{r_{1,1}}\dotsm x_{k,i}^{r_{k,i}+1}x_{k,i+1}^{r_{k,i+1}-1}\dotsm x_{m,n}^{r_{m,n}}
\end{align*}
and
\begin{align*}
  \tau_q(k_i)(x^r)&= q^{\sum_{t=1}^m (r_{t,i+1}-r_{t,i})}x^r
\end{align*}
for all $r \in M_{m,n}(\N_0)$ and $i=1,2,\dots,n-1$, which gives us \eqref{eq:rho action U_q(sl(n))}. Analogously, from \eqref{eq:adjoin action levi U_q(l)_b} and \eqref{eq:Hopf algebra structure U_q(l) b} we obtain
\begin{align*}
  \tau_q(e_{n+i})(x^r) &= \sum_{k=1}^n \sum_{j=1}^{r_{i+1,k}} q^{-\sum_{t=k+1}^n r_{i+1,k}-j+1} \,x_{1,1}^{r_{1,1}} \dotsm x_{i+1,k}^{r_{i+1,k}-j} x_{i,k} x_{i+1,k}^{j-1} \dotsm x_{m,n}^{r_{m,n}} \\
  &= \sum_{k=1}^n \sum_{j=1}^{r_{i+1,k}} q^{\sum_{t=k+1}^n(r_{i,k}-r_{i+1,k})+ r_{i+1,k}-2j+1} \,x_{1,1}^{r_{1,1}} \dotsm x_{i,k}^{r_{i,k}+1} \dotsm x_{i+1,k}^{r_{i+1,k}-1} \dotsm x_{m,n}^{r_{m,n}} \\
  &= \sum_{k=1}^n q^{\sum_{t=k+1}^n(r_{i,k}-r_{i+1,k})} [r_{i+1,k}]_q \,x_{1,1}^{r_{1,1}} \dotsm x_{i,k}^{r_{i,k}+1} \dotsm x_{i+1,k}^{r_{i+1,k}-1} \dotsm x_{m,n}^{r_{m,n}},
\end{align*}
\begin{align*}
  \tau_q(f_{n+i})(x^r) &= \sum_{k=1}^n \sum_{j=1}^{r_{i,k}} q^{-\sum_{t=1}^{k-1} r_{i,t}-j+1} \,x_{1,1}^{r_{1,1}} \dotsm x_{i,k}^{j-1} x_{i+1,k} x_{i,k}^{r_{i,k}-j} \dotsm x_{m,n}^{r_{m,n}} \\
  & =\sum_{k=1}^n \sum_{j=1}^{r_{i,k}} q^{\sum_{t=1}^{k-1} (r_{i+1,t}-r_{i,t})+ r_{i,k}-2j+1} \,x_{1,1}^{r_{1,1}} \dotsm x_{i,k}^{r_{i,k}-1} \dotsm x_{i+1,k}^{r_{i+1,k}+1} \dotsm x_{m,n}^{r_{m,n}} \\
  & =\sum_{k=1}^n q^{\sum_{t=1}^{k-1} (r_{i+1,t}-r_{i,t})} [r_{i,k}]_q \,x_{1,1}^{r_{1,1}} \dotsm x_{i,k}^{r_{i,k}-1} \dotsm x_{i+1,k}^{r_{i+1,k}+1} \dotsm x_{m,n}^{r_{m,n}}
\end{align*}
and
\begin{align*}
  \tau_q(k_{n+i})(x^r) &= q^{\sum_{t=1}^n (r_{i,t}-r_{i+1,t})} x^r
\end{align*}
for all $r \in M_{m,n}(\N_0)$ and $i=1,2,\dots,m-1$, which implies \eqref{eq:rho action U_q(sl(m))}. Finally, using \eqref{eq:adjoin action levi U_q(l)_c} and \eqref{eq:Hopf algebra structure U_q(l) c} we get
\begin{align*}
  \tau_q(k_n)(x^r) = q^{-\sum_{t=1}^n r_{1,t} -\sum_{s=1}^m r_{s,n}} x^r
\end{align*}
for all $r \in M_{m,n}(\N_0)$, which finishes the proof.}

\lemma{\label{lem:rho action}
We have
\begin{align}
\rho_q(E_{j,i})= - \sum_{1 \leq k_1 \leq \dots \leq k_s \leq m} (q^{-1}-q)^{\tau(k_1,\dots,k_s)-1} x_{k_1,i} \theta_{k_1,\dots,k_s} \partial_{k_s,j} \prod_{a=1}^s \prod_{t=1}^{k_a} \gamma_{t,i+a} \gamma_{t,i+a-1}^{-1}
\end{align}
for $1 \leq i<j\leq n$, where $s=j-i$, $\tau(k_1,\dots,k_s)$ is the number of distinct integers in the $s$-tuple $(k_1,k_2,\dots,k_s)$,
\begin{align}
  \theta_{k_1,\dots,k_s} = \beta_{k_1}\beta_{k_2} \dotsm \beta_{k_{s-1}}
\end{align}
with $\beta_{k_t}=\partial_{k_t,i+t} x_{k_{t+1},i+t}$ if $k_t \neq k_{t+1}$ and $\beta_{k_t}=\gamma_{k_t,i+t}^{-1}$ if $k_t= k_{t+1}$.}

\proof{We prove the statement by induction on $j-i$. The case $j-i=1$ follows from Theorem \ref{thm:rho action}. Further, for $s=j-i>1$ we have $\rho_q(E_{j,i})=\rho_q(E_{j,j-1}) \rho_q(E_{j-1,i})-q^{-1}\rho_q(E_{j-1,i})\rho_q(E_{j,j-1})$. By induction assumption we have
\begin{align*}
\rho_q(E_{j,j-1}) = -\sum_{k=1}^{m}x_{k,j-1}\p_{k,j}\prod_{t=1}^k\g_{t,j}\g_{t,j-1}^{-1}
\end{align*}
and
\begin{align*}
\rho_q(E_{j-1,i}) =  -\sum_{1 \leq k_1 \leq \dots \leq k_{s-1} \leq m}  (q^{-1}-q)^{\tau'-1}
 x_{k_1,i} \theta_{k_1,\dots,k_{s-1}} \partial_{k_{s-1},j-1} \prod_{a=1}^{s-1} \prod_{t=1}^{k_a} \gamma_{t,i+a} \gamma_{t,i+a-1}^{-1},
\end{align*}
where $\tau'=\tau(k_1,\dots,k_{s-1})$ for greater clarity, which gives us
\begin{align*}
 \rho_q(E_{j,i}) = \sum_{k=1}^m \sum_{1 \leq k_1 \leq \dots \leq k_{s-1} \leq m} (q^{-1}-q)^{\tau(k_1,\dots,k_{s-1})-1} \rho_q(E_{j,i})_{k,k_1,\dots,k_{s-1}},
\end{align*}
where $\rho_q(E_{j,i})_{k,k_1,\dots,k_{s-1}}$ denotes the expression
\begin{multline*}
  x_{k,j-1} \partial_{k,j} \prod_{t=1}^k \gamma_{t,j} \gamma_{t,j-1}^{-1}  x_{k_1,i} \theta_{k_1,\dots,k_{s-1}} \partial_{k_{s-1},j-1} \prod_{a=1}^{s-1} \prod_{t=1}^{k_a} \gamma_{t,i+a} \gamma_{t,i+a-1}^{-1} \\
 -q^{-1} x_{k_1,i} \theta_{k_1,\dots,k_{s-1}} \partial_{k_{s-1},j-1} \prod_{a=1}^{s-1} \prod_{t=1}^{k_a} \gamma_{t,i+a} \gamma_{t,i+a-1}^{-1} x_{k,j-1} \partial_{k,j} \prod_{t=1}^k \gamma_{t,j} \gamma_{t,j-1}^{-1}
\end{multline*}
for $1\leq k_1 \leq \dots \leq k_{s-1} \leq m$ and $k=1,2,\dots,m$.

If $k<k_{s-1}$, we have $\rho_q(E_{j,i})_{k,k_1,\dots,k_{s-1}} = 0$. For $k=k_{s-1}$ we may write
\begin{align*}
 \rho_q(E_{j,i})_{k,k_1,\dots,k_{s-1}} &=
  q  x_{k_1,i} \theta_{k_1,\dots,k_{s-1}} x_{k_{s-1},j-1} \partial_{k_{s-1},j-1} \partial_{k_{s-1},j} \prod_{a=1}^{s-1} \prod_{t=1}^{k_a} \gamma_{t,i+a} \gamma_{t,i+a-1}^{-1} \prod_{t=1}^{k_{s-1}} \gamma_{t,j} \gamma_{t,j-1}^{-1} \\
& \quad -  x_{k_1,i} \theta_{k_1,\dots,k_{s-1}} \partial_{k_{s-1},j-1} x_{k_{s-1},j-1} \partial_{k_{s-1},j} \prod_{a=1}^{s-1} \prod_{t=1}^{k_a} \gamma_{t,i+a} \gamma_{t,i+a-1}^{-1}
\prod_{t=1}^{k_{s-1}} \gamma_{t,j} \gamma_{t,j-1}^{-1} \\
&= - x_{k_1,i} \theta_{k_1,\dots,k_{s-1}} \gamma_{k_{s-1},j-1}^{-1} \partial_{k_{s-1},j} \prod_{a=1}^{s-1} \prod_{t=1}^{k_a} \gamma_{t,i+a} \gamma_{t,i+a-1}^{-1} \prod_{t=1}^{k_{s-1}} \gamma_{t,j} \gamma_{t,j-1}^{-1} \\
&= - x_{k_1,i} \theta_{k_1,\dots,k_{s-1},k_{s-1}} \partial_{k_{s-1},j} \prod_{a=1}^{s-1} \prod_{t=1}^{k_a} \gamma_{t,i+a} \gamma_{t,i+a-1}^{-1} \prod_{t=1}^{k_{s-1}} \gamma_{t,j} \gamma_{t,j-1}^{-1},
\end{align*}
where the second equality follows from $\partial_{i,j} x_{i,j} - q x_{i,j} \partial_{i,j} = \gamma_{i,j}^{-1}$. Finally, if $k>k_{s-1}$, we obtain
\begin{align*}
 \rho_q(E_{j,i})_{k,k_1,\dots,k_{s-1}}\! &= q x_{k_1,i} \theta_{k_1,\dots,k_{s-1}} \partial_{k_{s-1},j-1} x_{k,j-1}\partial_{k,j} \prod_{a=1}^{s-1} \prod_{t=1}^{k_a} \gamma_{t,i+a} \gamma_{t,i+a-1}^{-1} \prod_{t=1}^k \gamma_{t,j} \gamma_{t,j-1}^{-1} \\
& \quad -q^{-1} x_{k_1,i} \theta_{k_1,\dots,k_{s-1}} \partial_{k_{s-1},j-1} x_{k,j-1} \partial_{k,j}
 \prod_{a=1}^{s-1} \prod_{t=1}^{k_a} \gamma_{t,i+a} \gamma_{t,i+a-1}^{-1} \prod_{t=1}^k \gamma_{t,j} \gamma_{t,j-1}^{-1} \\
& =(q-q^{-1}) x_{k_1,i} \theta_{k_1,\dots,k_{s-1}} \partial_{k_{s-1},j-1} x_{k,j-1} \partial_{k,j}  \prod_{a=1}^{s-1} \prod_{t=1}^{k_a} \gamma_{t,i+a} \gamma_{t,i+a-1}^{-1} \prod_{t=1}^k \gamma_{t,j} \gamma_{t,j-1}^{-1} \\
&= (q-q^{-1}) x_{k_1,i} \theta_{k_1,\dots,k_{s-1},k}  \partial_{k,j}  \prod_{a=1}^{s-1} \prod_{t=1}^{k_a} \gamma_{t,i+a} \gamma_{t,i+a-1}^{-1} \prod_{t=1}^k \gamma_{t,j} \gamma_{t,j-1}^{-1}.
\end{align*}
Therefore, we have
\begin{align*}
  \rho_q(E_{j,i}) &= \sum_{k=k_{s-1}}^m \sum_{1 \leq k_1 \leq \dots \leq k_{s-1} \leq m} (q^{-1}-q)^{\tau(k_1,\dots,k_{s-1})-1} \rho_q(E_{j,i})_{k,k_1,\dots,k_{s-1}} \\
  &= - \sum_{1 \leq k_1 \leq \dots \leq k_s \leq m} (q^{-1}-q)^{\tau(k_1,\dots,k_s)-1} x_{k_1,i} \theta_{k_1,\dots,k_s} \partial_{k_s,j} \prod_{a=1}^s \prod_{t=1}^{k_a} \gamma_{t,i+a} \gamma_{t,i+a-1}^{-1},
\end{align*}
which gives the required statement.}

Now, let $(\sigma_q,V)$ be a $U_q(\mfrak{p})$-module. Then we can identify $U_q(\widebar{\mathfrak{u}})\otimes_\C\! V$ with $\C_q[\widebar{\mfrak{u}}^*] \otimes_\C\! V$ and obtain a $U_q(\mfrak{g})$-module structure on $\C_q[\widebar{\mfrak{u}}^*] \otimes_\C\! V$. Further, using the isomorphism $\varphi_q \colon \C[\widebar{\mfrak{u}}^*] \rarr \C_q[\widebar{\mfrak{u}}^*]$ of vector spaces, we can transfer the $U_q(\mfrak{g})$-module structure even on $\C[\widebar{\mfrak{u}}^*] \otimes_\C \! V$.

The main result of the present article is an explicit realization of the induced $U_q(\mfrak{g})$-module structure on $\C[\widebar{\mfrak{u}}^*] \otimes_\C \! V$ using quantum differential operators through the homomorphism
\begin{align}
  \pi_{q,V} \colon U_q(\mfrak{g}) \rarr \eus{A}^q_{\widebar{\mfrak{u}}^*}\! \otimes_\C \End V
\end{align}
of $\C$-algebras defined by
\begin{align}
  ((\psi_q \circ \varphi_q) \otimes \id_V)(\pi_{q,V}(a)(x^r\! \otimes v)) = a(E^r \otimes v) \label{eq:pi action definition}
\end{align}
for all $r \in M_{m,n}(\N_0)$ and $v \in V$. This is the content of the following theorem. Let us recall that $\C[\widebar{\mfrak{u}}^*]$ has a canonical structure of an $\eus{A}^q_{\widebar{\mfrak{u}}^*}$-module.

Since the $\C$-algebra isomorphism $\psi_q \colon \C_q[\widebar{\mfrak{u}}^*] \rarr U_q(\widebar{\mfrak{u}})$ is a homomorphism of $U_q(\mfrak{l})$-modules by construction of $\C_q[\widebar{\mfrak{u}}^*]$ and the vector space isomorphism $\varphi_q \colon \C[\widebar{\mfrak{u}}^*] \rarr \C_q[\widebar{\mfrak{u}}^*]$ is a homomorphism of $U_q(\mfrak{l})$-modules by definition of the $U_q(\mfrak{l})$-module structure on $\C[\widebar{\mfrak{u}}^*]$, the identification of $U_q(\widebar{\mfrak{u}}) \otimes_\C \! V$ with $\C_q[\widebar{\mfrak{u}}^*] \otimes_\C \! V$ is also an isomorphism of $U_q(\mfrak{l})$-modules and we obtain
\begin{align}
  \pi_{q,V}(a) = \sum \rho_q(a_{(1)}) \otimes \sigma_q(a_{(2)}),
\end{align}
where
\begin{align}
  \Delta(a) = \sum a_{(1)} \otimes a_{(2)},
\end{align}
for all $a \in U_q(\mfrak{l})$. If $V$ is the trivial $U_q({\mathfrak{p}})$-module, then we have $\pi_{q,V}(a) = \rho_q(a)$ for all $a \in U_q(\mfrak{l})$. In that case, we shall denote $\pi_{q,V}$ by $\rho_q$.
\medskip

\theorem{\label{thm}
Let $(\sigma_q,V)$ be a $U_q({\mathfrak{p}})$-module. Then the induced $U_q(\mathfrak{sl}_{n+m}(\C))$-module structure on $\C[\widebar{\mfrak{u}}^*] \otimes_\C \!V$ is defined through the homomorphism
\begin{align}
\pi_{q,V} \colon U_q(\mathfrak{sl}_{n+m}(\C)) \rarr \eus{A}^q_{\widebar{\mfrak{u}}^*} \!\otimes_\C \End V
\end{align}
of $\C$-algebras by
\begin{align}
\begin{aligned}
 \pi_{q,V}(f_i) &= \rho_q(f_i) \otimes \id_V + \rho_q(k_i) \otimes \sigma_q(f_i), \\
 \pi_{q,V}(e_i) &= \rho_q(e_i) \otimes \sigma_q(k_i^{-1}) + 1 \otimes \sigma_q(e_i), \\
 \pi_{q,V}(k_i) &= \rho_q(k_i) \otimes \sigma_q(k_i)
\end{aligned}
\end{align}
for $i=1,2,\dots,n-1$,
\begin{align}
\begin{aligned}
 \pi_{q,V}(f_{n+i}) &= \rho_q(f_{n+i}) \otimes \id_V + \rho_q(k_{n+i}^{-1}) \otimes \sigma_q(f_{n+i}), \\
 \pi_{q,V}(e_{n+i}) &= \rho_q(e_{n+i}) \otimes \sigma_q(k_{n+i}) + 1 \otimes \sigma_q(e_{n+i}), \\
 \pi_{q,V}(k_{n+i}) &= \rho_q(k_{n+i}) \otimes \sigma_q(k_{n+i})
\end{aligned}
\end{align}
for $i=1,2,\dots,m-1$, and
\begin{align}
\begin{aligned}
 \pi_{q,V}(f_n) &= x_{1,n} \prod_{t=1}^{n-1} \gamma_{1,t} \otimes \id_V, \\
 \pi_{q,V}(e_n) &= \sum_{k=1}^{n-1} \pi_{q,V}(E_{n,k}k_n)  \bigg(\prod_{t=1}^k \gamma_{1,t} \partial_{1,k} \otimes \id_V\!\bigg) - \sum_{k=2}^m \prod_{t=k}^m \gamma_{t,n}  \partial_{k,n} \otimes \sigma_q(k_n^{-1}E_{n+k,n+1}) \\
 & \quad - \sum_{k=1}^m  \prod_{t=1}^{k-1} \gamma_{t,n} \prod_{t=k+1}^m \gamma_{t,n}^{-1}\,  x_{k,n} \partial_{k,n} \partial_{1,n} \otimes \sigma_q(k_n) +  1  \otimes \sigma_q(e_n) \\
 & \quad + \prod_{t=1}^m \gamma_{t,n} \partial_{1,n} \otimes {\sigma_q(k_n) - \sigma_q(k_n^{-1}) \over q - q^{-1}}  \\
 \pi_{q,V}(k_n) &= \rho_q(k_n) \otimes \sigma_q(k_n),
\end{aligned}
\end{align}
where
\begin{align*}
\pi_{q,V}(E_{n,k})  = \rho_q(E_{n,k})\otimes \id_V + \rho_q(K_{k,n}) \otimes \sigma_q(E_{n,k}) + (q-q^{-1}) \sum_{k<\ell<n} \rho_q(E_{\ell,k}K_{\ell,n}) \otimes \sigma_q(E_{n,\ell})
\end{align*}
for $k=1,2,\dots,n-1$.}

\proof{From the previous considerations we know that the action of the Levi quantum subgroup $U_q(\mfrak{l})$ on $\C[\widebar{\mfrak{u}}^*] \otimes_\C \! V$ is given through the homomorphism
\begin{align*}
  \pi_{q,V} \colon U_q({\mfrak{l}}) \rarr \eus{A}^q_{\widebar{\mfrak{u}}^*} \otimes \End V
\end{align*}
of $\C$-algebras by the formula
\begin{align*}
  \pi_{q,V}(a) = \sum \rho_q(a_{(1)}) \otimes \sigma_q(a_{(2)}),
\end{align*}
where $\Delta(a) = \sum a_{(1)} \otimes a_{(2)}$, for all $a \in U_q(\mfrak{l})$. Hence, by using Theorem \ref{thm:rho action} we get the corresponding expressions for all generators of $U_q(\mfrak{g})$ except $e_n$ and $f_n$.

By Lemma \ref{lem:relations} we have
\begin{align*}
  E_{n+1,n}E^r = q^{\sum_{t=1}^{n-1} \!r_{1,t}} E^{r+1_{1,n}},
\end{align*}
which together with \eqref{eq:pi action definition} gives us
\begin{align*}
  \pi_{q,V}(E_{n+1,n})(x^r\! \otimes v) = q^{\sum_{t=1}^{n-1} \!r_{1,t}} x^{r+1_{1,n}} \otimes v  = \bigg(x_{1,n} \prod_{t=1}^{n-1} \gamma_{1,t} \otimes \id_V\!\bigg) (x^r\! \otimes v)
\end{align*}
for all $r \in M_{m,n}(\N_0)$ and $v \in V$. Further, using Lemma \ref{lem:relations} we may write
\begin{align*}
  [E_{n,n+1},E_{n+1,1}^{r_{1,1}}\dotsm E_{n+m,n}^{r_{m,n}}] &= E_{n+1,1}^{r_{1,1}} \dotsm [E_{n,n+1}, E_{n+1,n}^{r_{1,n}}] \dotsm E_{n+m,n}^{r_{m,n}} \\
  & \quad + \sum_{k=1}^{n-1} E_{n+1,1}^{r_{1,1}} \dotsm [E_{n,n+1}, E_{n+1,k}^{r_{1,k}}] \dotsm E_{n+m,n}^{r_{m,n}} \\
  & \quad + \sum_{k=2}^m E_{n+1,1}^{r_{1,1}} \dotsm [E_{n,n+1}, E_{n+k,n}^{r_{k,n}}] \dotsm E_{n+m,n}^{r_{m,n}} \\
  &=  [r_{1,n}]_q E^{r-1_{1,n}}  {q^{-\sum_{t=1}^m\! r_{t,n} +1}K_{n,n+1} -  q^{\sum_{t=1}^m \! r_{t,n} -1} K_{n,n+1}^{-1} \over q - q^{-1}} \\
  & \quad + \sum_{k=1}^{n-1}\, [r_{1,k}]_q q^{-\sum_{t=k+1}^n \!r_{1,t} - \sum_{t=1}^m \!r_{t,n}} E_{n,k} E^{r-1_{1,k}} K_{n,n+1} \\
  & \quad - \sum_{k=2}^m q^{-2} [r_{k,n}]_q q^{\sum_{t=k}^m\! r_{t,n}} E^{r-1_{k,n}} E_{n+k,n+1} K_{n,n+1}^{-1},
\end{align*}
which together with \eqref{eq:pi action definition} implies
\begin{align*}
     \pi_{q,V}(E_{n,n+1})(x^r\! \otimes v) & = x^r \!\otimes \sigma_q(E_{n,n+1})v  \\
     & \quad + [r_{1,n}]_q x^{r-1_{1,n}} \! \otimes {q^{-\sum_{t=1}^m\! r_{t,n} +1}\,\sigma_q(K_{n,n+1}) -  q^{\sum_{t=1}^m \! r_{t,n} -1}\, \sigma_q(K_{n,n+1}^{-1}) \over q - q^{-1}}\, v \\
  & \quad + \sum_{k=1}^{n-1}\, [r_{1,k}]_q q^{-\sum_{t=k+1}^n \!r_{1,t} - \sum_{t=1}^m \!r_{t,n}}\, \pi_{q,V}(E_{n,k})(x^{r-1_{1,k}}\! \otimes \sigma_q(K_{n,n+1})v) \\
  & \quad - \sum_{k=2}^m q^{-2} [r_{k,n}]_q q^{\sum_{t=k}^m\! r_{t,n}} x^{r-1_{k,n}}\! \otimes \sigma_q(E_{n+k,n+1} K_{n,n+1}^{-1})v
\end{align*}
for all $r \in M_{m,n}(\N_0)$ and $v \in V$. Furthermore, using
\begin{align*}
    {q^{-\sum_{t=1}^m\! r_{t,n} +1}K_{n,n+1} -  q^{\sum_{t=1}^m \! r_{t,n} -1} K_{n,n+1}^{-1} \over q - q^{-1}} & = q^{\sum_{t=1}^m \! r_{t,n}-1} {K_{n,n+1} - K_{n,n+1}^{-1} \over q - q^{-1}} \\
    & \quad -  {q^{\sum_{t=1}^m \! r_{t,n}-1} -  q^{-\sum_{t=1}^m \!r_{t,n}+1} \over q - q^{-1}}\, K_{n,n+1},
\end{align*}
we get
\begin{align*}
  \pi_{q,V}(E_{n,n+1})(x^r\! \otimes v) & = x^r \!\otimes \sigma_q(E_{n,n+1})v + \prod_{t=1}^m \gamma_{t,n} \partial_{1,n}x^r \! \otimes {\sigma_q(K_{n,n+1}) - \sigma_q(K_{n,n+1}^{-1}) \over q - q^{-1}}\, v  \\
  & \quad - {\prod_{t=1}^m \gamma_{t,n} - \prod_{t=1}^m \gamma_{t,n}^{-1} \over q - q^{-1}}\, \partial_{1,n} x^r \! \otimes \sigma_q(K_{n,n+1})v \\
  & \quad + \sum_{k=1}^{n-1} \pi_{q,V}(E_{n,k})\bigg(\rho_q(K_{n,n+1}) \prod_{t=1}^k \gamma_{1,t} \partial_{1,k} x^r \! \otimes \sigma_q(K_{n,n+1})v\!\bigg) \\
  & \quad - \sum_{k=2}^m  \prod_{t=k}^m \gamma_{t,n} \partial_{k,n} x^r \! \otimes \sigma_q(K_{n,n+1}^{-1} E_{n+k,n+1})v
\end{align*}
for all $r \in M_{m,n}(\N_0)$ and $v \in V$. Finally, using the formula
\begin{align*}
  {\prod_{t=1}^m \gamma_{t,n} - \prod_{t=1}^m \gamma_{t,n}^{-1} \over q - q^{-1}} = \sum_{k=1}^m  \prod_{t=1}^{k-1} \gamma_{t,n} \prod_{t=k+1}^m \gamma_{t,n}^{-1}  x_{k,n} \partial_{k,n},
\end{align*}
we obtain the required statement. Moreover, from Lemma \ref{lem:coproduct} we have
\begin{align*}
\pi_{q,V}(E_{n,k}) = \rho_q(E_{n,k})\otimes \id_V + \rho_q(K_{k,n}) \otimes \sigma_q(E_{n,k})
 + (q-q^{-1}) \sum_{k<\ell<n} \rho_q(E_{\ell,k}K_{\ell,n}) \otimes \sigma_q(E_{n,\ell})
\end{align*}
for $k=1,2,\dots,n-1$. This finishes the proof.}

Theorem \ref{thm} together with Theorem \ref{thm:rho action} and Lemma \ref{lem:rho action} give us an explicit realization of induced modules $M^\mfrak{g}_{q,\mfrak{p}}(V) \simeq \C[\widebar{\mfrak{u}}^*] \otimes_\C V$ by quantum differential operators for any $U_q(\mfrak{p})$-module $V$. In fact, we have a stronger result. As $\pi_{q,V} \colon U_q(\mfrak{g}) \rarr \eus{A}^q_{\widebar{\mfrak{u}}^*} \!\otimes_\C \End V$ is a homomorphism of associative $\C$-algebras, we may take another $\eus{A}^q_{\widebar{\mfrak{u}}^*}$-module $\eus{M}$ instead of $\C[\widebar{\mfrak{u}}^*]$ and we obtain a $U_q(\mfrak{g})$-module structure on $\eus{M} \otimes_\C\! V$ through the homomorphism $\pi_{q,V}$. Moreover, since the classical limit of $\eus{A}^q_{\widebar{\mfrak{u}}^*}$ via the specialization $q \rarr 1$ is the Weyl algebra $\eus{A}_{\widebar{\mfrak{u}}^*}$, it would be interesting to consider such $\eus{A}^q_{\widebar{\mfrak{u}}^*}$-modules for $\eus{M}$ that the corresponding $U_q(\mfrak{g})$-module $\eus{M} \otimes_\C \! V$ is a true deformation of a twisted induced module.
\medskip

A finite-dimensional module $V$ over $U_q(\mathfrak{l})$ has the classical limit $\widetilde{V}$ over $\mathfrak{l}$. As it was mentioned at the end of the previous section, extending properly $V$ to a module over $U_q(\mathfrak{p})$ we can guarantee that it still admits the classical limit and that $\C_q[\widebar{\mfrak{u}}^*] \otimes_\C\! V$ is a true deformation of a generalized Verma module for $\mathfrak{g}$. Moreover, from Theorem \ref{thm} we easily get the classical limit by the specialization $q\rarr 1$.


\section*{Acknowledgments}

V.\,F.\ is supported in part by CNPq (304467/2017-0) and by Fapesp (2014/09310-5); L.\,K.\ is supported by Capes (88887.137839/2017-00) and J.\,Z.\ is supported by Fapesp (2015/05927-0).



\begin{thebibliography}{DMP00}

\bibitem[BC90]{Boe-Collingwood1990}
Brian Boe and David Collingwood, \emph{{Multiplicity free categories of highest
  weight representations. I}}, Comm. Algebra \textbf{18} (1990), no.~4,
  947--1032.

\bibitem[CF94a]{Coleman-Futorny1994}
John~A. Coleman and Vyacheslav Futorny, \emph{{Stratified $L$-modules}}, J.
  Algebra \textbf{163} (1994), no.~1, 219--234.

\bibitem[CP94b]{Chari-Pressley1994}
Vyjayanthi Chari and Andrew Pressley, \emph{{A guide to quantum groups}},
  Cambridge University Press, Cambridge, 1994.

\bibitem[DMP00]{Dimitrov-Mathieu-Penkov2000}
Ivan Dimitrov, Olivier Mathieu, and Ivan Penkov, \emph{{On the structure of
  weight modules}}, Trans. Amer. Math. Soc. \textbf{352} (2000), no.~6,
  2857--2869.

\bibitem[DOF90]{Drozd-Ovsienko-Futorny1990}
Yuri~A. Drozd, Serge Ovsienko, and Vyacheslav Futorny, \emph{{The
  Harish-Chandra $S$-homomorphism and $\mathfrak{G}$-modules generated by
  primitive elements}}, Ukrainian Math. J. \textbf{42} (1990), no.~8, 919--924.

\bibitem[Fer90]{Fernando1990}
Suren~L. Fernando, \emph{{Lie algebra modules with finite-dimensional weight
  spaces. I}}, Trans. Amer. Math. Soc. \textbf{322} (1990), no.~2, 757--781.

\bibitem[Fut86]{Futorny1986}
Vyacheslav Futorny, \emph{{A generalization of Verma modules and irreducible
  representations of the Lie algebra $\mathfrak{sl}(3)$}}, Ukrain. Mat. Zh.
  \textbf{38} (1986), no.~4, 492--497.

\bibitem[Fut87]{Futorny1987}
\bysame, \emph{{The weight representations of semisimple finite dimensional Lie
  algebras}}, Ph.D. thesis, Kiev University, 1987.

\bibitem[GL76]{Garland-Lepowsky1976}
Howard Garland and James Lepowsky, \emph{{Lie algebra homology and the
  Macdonald-Kac formulas}}, Invent. Math. \textbf{34} (1976), no.~1, 37--76.

\bibitem[IS88]{Irving-Shelton1988}
Ronald~S. Irving and Brad Shelton, \emph{{Loewy series and simple projective
  modules in the category $\mathcal{O}_S$}}, Pacific J. Math. \textbf{132}
  (1988), no.~2, 319--342.

\bibitem[Jim86]{Jimbo1986}
Michio Jimbo, \emph{{A $q$-analogue of $U(\mathfrak{gl}(N+1))$, Hecke algebra,
  and the Yang--Baxter equation}}, Lett. Math. Phys. \textbf{11} (1986), no.~3,
  247--252.

\bibitem[Kas95]{Kassel1995}
Christian Kassel, \emph{{Qantum groups}}, Graduate Texts in Mathematics, vol.
  155, Springer-Verlag, New York, 1995.

\bibitem[KS97]{Klimyk-Schmudgen1997}
Anatoli Klimyk and Konrad Schmüdgen, \emph{{Qantum groups and their
  representations}}, Springer-Verlag, Berlin, 1997.

\bibitem[Lep77a]{Lepowsky1977}
James Lepowsky, \emph{{A generalization of the Bernstein-Gelfand-Gelfand
  resolution}}, J. Algebra \textbf{49} (1977), no.~2, 496--511.

\bibitem[Lep77b]{Lepowsky1977a}
\bysame, \emph{{Generalized Verma modules, the Cartan-Helgason theorem, and the
  Harish-Chandra homomorphism}}, J. Algebra. \textbf{49} (1977), no.~2,
  470--495.

\bibitem[Lus88]{Lusztig1988}
George Lusztig, \emph{{Quantum deformations of certain simple modules over
  enveloping algebras}}, Adv. Math. \textbf{70} (1988), no.~2, 237--249.

\bibitem[Maz00]{Mazorchuk2000}
Volodymyr Mazorchuk, \emph{{Generalized Verma modules}}, Mathematical Studies
  Monograph Series, vol.~8, VNTL Publishers, Lviv, 2000.

\bibitem[Mel99]{Melville1999}
Duncan~J. Melville, \emph{{An $\mathbb{A}$-form technique of quantum
  deformations}}, Recent developments in quantum affine algebras and related
  topics (Raleigh, NC, 1998), Contemporary Mathematics, vol. 248, Amer. Math.
  Soc., Providence, RI, 1999, pp.~359--375.

\bibitem[MS08]{Mazorchuk-Stroppel2008}
Volodymyr Mazorchuk and Catharina Stroppel, \emph{{Categorification of
  (induced) cell modules and the rough structure of generalised Verma
  modules}}, Adv. Math. \textbf{219} (2008), no.~4, 1363--1426.

\bibitem[RC80]{Rocha-Caridi1980}
Alvany Rocha-Caridi, \emph{{Splitting criteria for $\mathfrak{g}$-modules
  induced from a parabolic and the Bernstein-Gelfand-Gelfand resolution of a
  finite-dimensional, irreducible $\mathfrak{g}$-modules}}, Trans. Amer. Math.
  Soc. \textbf{262} (1980), no.~2, 335--366.

\bibitem[RTF90]{Reshetikhin-Takhtajan-Faddeev1990}
Nikolai~Yu. Reshetikhin, Leon~A. Takhtajan, and Lyudvig~D. Faddeev,
  \emph{{Quantization of Lie groups and Lie algebras}}, Leningrad Math. J.
  \textbf{1} (1990), no.~1, 193--225.
  
\bibitem[S86]{Shen1986}  G. Shen,   \emph{Graded modules of graded Lie algebras of Cartan type. I. Mixed
products of modules}, Sci. Sinica Ser., \textbf{A 29} (1986), no. 6, 570-581.

\bibitem[Ver66]{Verma1966}
Dayad-Nand Verma, \emph{{Structure of certain induced representations of
  complex semisimple Lie algebras}}, Ph.D. thesis, Yale University, 1966.

\end{thebibliography}

\providecommand{\bysame}{\leavevmode\hbox to3em{\hrulefill}\thinspace}
\providecommand{\MR}{\relax\ifhmode\unskip\space\fi MR }
\providecommand{\MRhref}[2]{%
  \href{http://www.ams.org/mathscinet-getitem?mr=#1}{#2}
}
\providecommand{\href}[2]{#2}

\end{document}